\documentclass[11pt]{amsart}
\usepackage{amsthm,amsmath,amssymb, stmaryrd}
\usepackage[dvipsnames,svgnames,x11names,hyperref]{xcolor}
\usepackage[T1]{fontenc}
\usepackage[english]{babel}
\usepackage[numbers]{natbib}
\usepackage[colorlinks=true,citecolor=Cerulean]{hyperref}
\usepackage{graphicx}
\usepackage{lipsum}
\usepackage{tikz}
\usepackage{caption}
\usepackage{subcaption}
\usepackage{array}
\newcolumntype{P}[1]{>{\centering\arraybackslash}p{#1}}
\usetikzlibrary{calc}
\newtheorem{theorem}{Theorem}[section]
\newtheorem{lemma}{Lemma}[section]
\newtheorem{definition}{Definition}[section]
\newtheorem{corollary}{Corollary}[section]
\newtheorem{proposition}{Proposition}[section]

\newtheorem{remark}{Remark}[section]
\newtheorem{notation}{Notation}[section]
\newcommand{\Le}{\rotatebox[origin=c]{180}{$\Gamma$}}
    \DeclareFontFamily{U}{wncy}{}
    \DeclareFontShape{U}{wncy}{m}{n}{<->wncyr10}{}
    \DeclareSymbolFont{mcy}{U}{wncy}{m}{n}
    \DeclareMathSymbol{\Sh}{\mathord}{mcy}{"58} 

\title{Excluded minors of interval positroids that are paving matroids}
\author{Hyungju Park}
\date{}

\begin{document}
\begin{abstract}

We prove that every paving matroid that is an excluded minor of interval positroids can be reduced to one of three fundamental families of excluded minors of interval positroids by relaxing dependent hyperplanes. Using this result, we classify all non-positroid excluded minors of interval positroids that are paving matroids. Additionally, we provide a criterion that characterizes all excluded minors of interval positroids that are paving positroids.
    
\end{abstract}
\maketitle
\section{Introduction}

\emph{Positroids} are matroids introduced and studied by \citet{Pos} and, under the name ``base-sortable matroids'', by \citet{MR1857256}. The term ``positroid'' is used to refer to two slightly different objects in the literature. A positroid might indicate a matroid $M$ on the totally ordered set $[n]:=\{1,\ldots,n\}$ that is represented by a real matrix with $n$ columns and nonnegative maximal minors, by the correspondence associating each $i\in [n]$ with the $i$th column. A matroid that is isomorphic to a positroid, according to the previous description, is also called a positroid. To avoid confusion, the positroid in the former sense will be called an \emph{ordered positroid}. \citet{Pos} found various combinatorial objects that bijectively correspond to ordered positroids such as Grassmann necklaces, decorated permutations, $\Le$-diagrams, and move-equivalence classes of reduced plabic graphs.

\citet{MR3413866} observed that the class of positroids is closed under taking minors and direct sum. Since the class of positroids is minor-closed, it is natural to consider excluded minors for the class of positroids. \citet[Section~4]{MR2721556} showed that an ordered matroid is an ordered positroid if and only if it avoids a minor with bases $\{ij,jk,kl,il\}$, $\{ij,jk,kl,il, ik\}$, or $\{ij,jk,kl,il,jl\}$ on any subset $\{i,j,k,l\}\subseteq [n]$ with $i<j<k<l$. But, excluded minors for the class of (unordered) positroids are not fully characterized yet. \citet{MR1857256} characterized all excluded minors of positroids of rank three. \citet[Section~5]{bonin2023characterization} found several infinite families of excluded minors for the class of positroids.

\emph{Lattice path matroids} (LPMs) are a special minor-closed class of matroids associated with any pair of two noncrossing lattice paths that share the same initial and terminal points. \citet{MR2834184} showed that LPMs are positroids and \citet{MR2718679} characterized excluded minors of LPMs. \emph{Interval positroids}, introduced by \citet{MR3466433}, constitute a minor-closed subclass of positroids that include the class of LPMs. We take the first step in studying excluded minors of interval positroids by characterizing all cases where they are also paving matroids. In the class of paving matroids, positroids precisely correspond to matroids with a circular-arc collection of dependent hyperplanes. Our goal is to characterize all hypergraphs that are collections of all dependent hyperplanes of a paving excluded minor of interval positroids.

In Section \ref{sec:back}, we review some results on interval and circular-arc hypergraphs, positroids, and the hyperplane characterization of paving matroids. In Section \ref{sec:rel}, three families of ``relaxation-minimal'' paving excluded minors of interval positroids are classified. These families are relaxation-minimal in the sense that every other paving excluded minor of interval positroids can be reduced to one of these three families through relaxation of stressed hyperplanes. In Section \ref{sec:pos}, we give a criterion for characterizing paving excluded minors of interval positroids that are also positroids. In Section \ref{sec:four}, every paving, non-positroid excluded minors of interval positroids with four dependent hyperplanes is determined. In Section \ref{sec:five}, we classify every paving, non-positroid excluded minors of interval positroids with more than four dependent hyperplanes. In Section \ref{sec:sum}, we summarize the results and identify which of these excluded minors have been discussed in the recent paper by \citet{bonin2023characterization}.

\section{Background}\label{sec:back}

For notions and terms related to matroids that are not explained here, we refer to the book by \citet{MR2849819}.

\subsection{Interval hypergraphs and CA hypergraphs}

Let us introduce the following terminology.

\begin{definition}

Two sets $A, B$ overlap if their intersection is nonempty and incomparable ($A\not\subset B, B\not\subset A$).
    
\end{definition}

A \emph{hypergraph} $\mathcal{H}$ is a pair $(V,E)$ of the set of vertices $V$ and the set of edges $E$ which is a collection of subsets of $V$. A hypergraph is a \emph{clutter} if any distinct two edges are incomparable. A \emph{partial hypergraph} of the hypergraph $\mathcal{H}=(V,E)$ is some hypergraph $\mathcal{H}'=(V,E')$ where the set of edges $E'$ is contained in the original set of edges $E$. A \emph{subhypergraph} of a hypergraph is the hypergraph $\mathcal{H}$ obtained by removing some vertices of a partial hypergraph of $\mathcal{H}$. A \emph{partial clutter} is a partial hypergraph that is a clutter. A vertex $x\in V$ is \emph{isolated} if it is not contained in any edge. 

\begin{figure}
    \centering
    \subcaptionbox{Type $1$}{
    \begin{tikzpicture}
      \draw[thick] (30:0.5) circle [radius=1];
    \draw[thick] (150:0.5) circle [radius=1];
    \draw[thick] (270:0.5) circle [radius=1];
    \draw[very thick] (90:0.8) node[circle, fill=black, inner sep=1.5pt] {};
    \draw[very thick] (330:0.8) node[circle, fill=black, inner sep=1.5pt] {};
    \draw[very thick] (210:0.8) node[circle, fill=black, inner sep=1.5pt] {};
    \end{tikzpicture}
    }
    \subcaptionbox{Type $2$}{
    \begin{tikzpicture}
      \draw[thick] (30:0.5) circle [radius=1];
    \draw[thick] (150:0.5) circle [radius=1];
    \draw[thick] (270:0.5) circle [radius=1];
    \draw[very thick] (30:1) node[circle, fill=black, inner sep=1.5pt] {};
    \draw[very thick] (270:1) node[circle, fill=black, inner sep=1.5pt] {};
    \draw[very thick] (150:1) node[circle, fill=black, inner sep=1.5pt] {};
    \draw[very thick] (0,0) node[circle, fill=black, inner sep=1.5pt] {};
    \end{tikzpicture}
    }
    \subcaptionbox{Type $3$}{
    \begin{tikzpicture}
      \draw[thick] (45:0.5) circle [radius=1];
    \draw[thick] (135:0.5) circle [radius=1];
    \draw[thick] (225:0.5) circle [radius=1];
    \draw[thick] (315:0.5) circle [radius=1];
    \draw[very thick] (0:1) node[circle, fill=black, inner sep=1.5pt] {};
    \draw[very thick] (90:1) node[circle, fill=black, inner sep=1.5pt] {};
    \draw[very thick] (180:1) node[circle, fill=black, inner sep=1.5pt] {};
    \draw[very thick] (270:1) node[circle, fill=black, inner sep=1.5pt] {};
    \draw[very thick] (270:0) node[circle, fill=black, inner sep=1.5pt] {};
    \end{tikzpicture}
    }
    \subcaptionbox{Type $4$}{
    \begin{tikzpicture}
      \draw[thick] (30:1) circle [radius=0.6];
    \draw[thick] (150:1) circle [radius=0.6];
    \draw[thick] (270:1) circle [radius=0.6];
    \draw[thick] (270:0) circle [radius=1];
    \draw[very thick] (30:0.7) node[circle, fill=black, inner sep=1.5pt] {};
    \draw[very thick] (30:1.3) node[circle, fill=black, inner sep=1.5pt] {};
    \draw[very thick] (270:0.7) node[circle, fill=black, inner sep=1.5pt] {};
    \draw[very thick] (270:1.3) node[circle, fill=black, inner sep=1.5pt] {};
    \draw[very thick] (150:0.7) node[circle, fill=black, inner sep=1.5pt] {};
    \draw[very thick] (150:1.3) node[circle, fill=black, inner sep=1.5pt] {};
    \end{tikzpicture}
    }
    \subcaptionbox{Type $5$}{
    \begin{tikzpicture}
    \draw[thick] (30:1) circle [radius=0.4];
    \draw[thick] (60:1) circle [radius=0.4];
    \draw[thick] (90:1) circle [radius=0.4];
    \draw[thick] (120:1) circle [radius=0.4];
    \draw[thick] (150:1) circle [radius=0.4];
    \draw[very thick] (45:1) node[circle, fill=black, inner sep=1.5pt] {};
    \draw[very thick] (75:1) node[circle, fill=black, inner sep=1.5pt] {};
    \draw[very thick] (105:1) node[circle, fill=black, inner sep=1.5pt] {};
    \draw[very thick] (135:1) node[circle, fill=black, inner sep=1.5pt] {};
    \fill (180:1) circle (1pt);
    \fill (190:1) circle (1pt);
    \fill (200:1) circle (1pt);
    \fill (360:1) circle (1pt);
    \fill (350:1) circle (1pt);
    \fill (340:1) circle (1pt);

    \end{tikzpicture}
    }
    \caption{Forbidden partial clutters of interval clutters}
    \label{fig:forint}
\end{figure}
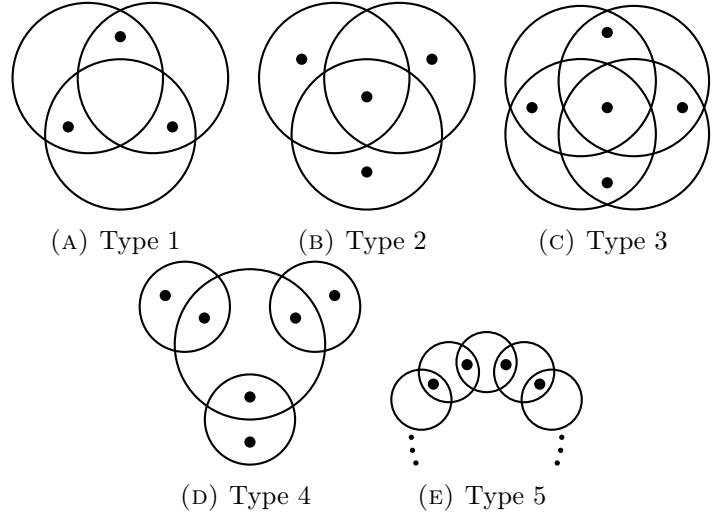

A hypergraph with $n$ vertices is an \emph{interval hypergraph} if there is a bijection from $V$ to the set $\{1,\ldots,n\}$ such that the image of every edge under this bijection is an interval $[i,j]$ for some $1\le i\le j\le n$. Such a bijection will be called an \emph{interval ordering} of the hypergraph. A hypergraph with $n$ vertices is a \emph{circular-arc} (or \emph{CA} for brevity) hypergraph if there is a bijection from $V$ to the set $\{1,\ldots,n\}$ such that the image of every edge under this bijection is a cyclic interval $[i,j]$ or $[1,i]\cup [j,n]$ for some $1\le i\le j\le n$. Such a bijection will be called an \emph{arc ordering} of the hypergraph. A partial hypergraph or a subhypergraph of an interval (resp. CA) hypergraph is also an interval (resp. CA) hypergraph. \citet[Theorem~2]{MR1909862} classified forbidden partial clutters for the class of interval clutters whose Venn diagrams are drawn in Figure \ref{fig:forint} with dots to indicate nonempty sets.

\begin{proposition}(\cite[Theorem~2]{MR1909862})\label{prop:int}

A clutter $\mathcal{H}=(V,E)$ is non-interval if and only if some partial clutter of $\mathcal{H}$ is one of the following five types of clutters:

\begin{enumerate}
    \item $A,B,C\in E$ such that $(A\cap B)\setminus C, (B\cap C)\setminus A$ and $(A\cap C)\setminus B$ are all nonempty.

    \item $A,B,C\in E$ such that $A\setminus (B\cup C), B\setminus (A\cup C), C\setminus(A\cup B)$ and $A\cap B\cap C$ are all nonempty.

    \item $A,B,C,D\in E$ such that $(A\cap B)\setminus (C\cup D), (B\cap C)\setminus (D\cup A), (C\cap D)\setminus (A\cup B), (D\cap A)\setminus (B\cup C)$ and $A\cap B\cap C\cap D$ are all nonempty.

    \item $A,B,C,D\in E$ such that $A\cap B, B\cap C$ and $A\cap C$ are empty and $A\cap D, B\cap D$ and $C\cap D$ are nonempty.

    \item For $m\ge 4$, $A_1,\ldots, A_m\in E$ such that $A_i\cap A_j\neq\emptyset$ if and only if $i=j$ or $i=j\pm 1 \mod m$.
\end{enumerate}
    
\end{proposition}

Throughout the paper, we will refer to clutters described in Proposition \ref{prop:int} as non-interval clutters of type 1, 2, 3, 4, and 5. Note that non-interval clutters of type 2,3, and 4 are not CA. Therefore, the following corollary is a direct consequence of Proposition \ref{prop:int}.

\begin{corollary}\label{cor:ca}

A non-interval, CA clutter $\mathcal{H}=(V,E)$ contains a partial non-interval clutter of type 1 or 5.

\end{corollary}

The following two lemmas on CA clutters will be used in next sections.

\begin{lemma}\label{lem:tri}

Let $\mathcal{H}$ be a clutter with four edges $A,B,C,D$ such that:

\begin{enumerate}
    \item $A\cap D$, $B\cap D$, and $C\cap D$ are all nonempty;

    \item $A\setminus D$, $B\setminus D$, and $C\setminus D$ are pairwise disjoint.
\end{enumerate} Then $\mathcal{H}$ is non-interval.

\begin{proof}

Suppose that there is an arc ordering of $\mathcal{H}$ such that $A,B,C,D$ are all cyclic intervals. Since $D$ is a cyclic interval and $A\cap D$ is nonempty, $A\setminus D$ is a nonempty cyclic interval adjacent to $D$. By the same reason, $B\setminus D$ and $C\setminus D$ are nonempty cyclic intervals adjacent to $D$. This is impossible since $A\setminus D$, $B\setminus D$, and $C\setminus D$ are pairwise disjoint.
    
\end{proof}
    
\end{lemma}

\begin{lemma}\label{lem:arcorder}

Let $\mathcal{H}=(V,E=\{A_1,\ldots,A_m\})$ be a non-interval CA clutter with every vertex contained in at most two edges. Then, every intersection $A_i\cap A_{i+1}$ and every set of isolated points $A_i\setminus (A_{i-1}\cup A_{i+1})$ map to cyclic intervals for all $i\in[m]$ by every arc ordering $f:V\to [n]$ of $\mathcal{H}$. Moreover, the cyclic intervals appear in the following order up to rotation and reversal: $f(A_1\cap A_m), f(A_1\setminus(A_m\cup A_2)), f(A_1\cap A_2), f(A_2\setminus (A_1\cup A_3)),\ldots, f(A_m\setminus (A_{m-1}\cup A_1))$.

\begin{proof}

For an arc ordering $f:V\to [n]$, the collection of all images $f(A_i)$ is a collection cyclic interval. Since $\mathcal{H}$ is not an interval clutter, $m\ge 3$. Therefore, for each $i\in[m]$, $f(A_{i-1}), f(A_i), f(A_{i+1})$ are distinct cyclic intervals of the set $[n]$. Since both $f(A_{i-1})\cap f(A_i)$ and $f(A_{i+1})\cap f(A_i)$ are nonempty and disjoint, they are cyclic intervals at each end of $f(A_i)$. It follows that $f(A_i\setminus(A_{i-1}\cup A_{i+1}))$ is the cyclic interval in the middle obtained by removing both $f(A_{i-1})\cap f(A_i)$ and $f(A_{i+1})\cap f(A_i)$ from $f(A_i)$.

\end{proof}
    
\end{lemma}

\subsection{Positroids}

A \emph{positroid} is a $\mathbb{R}$-representable matroid that is representable by a matrix whose maximal minors are all nonnegative. As mentioned in the introduction, a positroid on the set $[n]$, represented by a nonnegative maximal minor matrix with the correspondence $i\mapsto$ the $i$th column, will be called an \emph{ordered positroid}. \citet{MR3413866} observed that a matroid is an ordered positroid with the given order of $[n]$ if and only if every connected component is a positroid with the induced order and the collection of connected components forms a noncrossing partition of the set $[n]$. Therefore, every excluded minor of positroids is connected.

The matroid polytope of a matroid $M$ on the set $[n]$ is defined by the convex hull of all vectors $e_B\in\mathbb{R}^n$ where $e_B=\sum_{i\in B}e_i$. For a connected matroid, each facet of its matroid polytope containing some interior point of the hypersimplex $\Delta_{d,n}$ is called a \emph{split facet}. Each split facet of a connected matroid is determined by a dependent flat $F$ of $M$ and the inequality $\sum_{i\in F}x_i\le r_M(F)$, see \cite[Propsition~7]{MR3663497}. The flat that corresponds to a split facet is called the \emph{split flacet}. Positroids can be characterized in terms of matroid polytope. A matroid is an ordered positroid with the given order if and only if its matroid polytope is a subpolytope of the hypersimplex $\Delta_{d,n}$ that is determined by inequalities of the form $x_i+x_{i+1}+\cdots x_j\le k$, where $[i,j]$ is a cyclic interval of the set $[n]$, see \cite[Proposition~5.7]{MR3413866}. Thus, a connected matroid is an ordered positroid if and only if all split flacets are cyclic intervals. Therefore, a connected matroid $M$ is a positroid if and only if the collection of every split flacets of $M$ is a CA hypergraph.

An (ordered) \emph{interval positroid} is characterized by a matroid whose matroid polytope is the subpolytope of the hypersimplex $\Delta_{d,n}$ that is determined by inequalities of the form $x_i+x_{i+1}+\cdots x_j\le k$, where $[i,j]$ is an interval of the set $[n]$, introduced by \citet{MR3466433}. Any matroid that is isomorphic to an ordered interval positroid will be called an interval positroid. The class of interval positroids is closed under taking minors but not under taking the dual.

\subsection{Hypergraphic characterization of paving matroids}

A \emph{paving matroid} can be characterized in several different ways. For our purpose, the most suitable approach is to characterize all hypergraphs that can be the set of hyperplanes of a paving matroid. \citet{MR0099931} introduced the notion of $d$-partition, which generalizes the usual set partition. A hypergraph with the vertex set $[n]$ is the collection of hyperplanes of a paving matroid of rank $d$ if and only if it is a $(d-1)$-partition of the set $[n]$.

\begin{definition}(\cite[Proposition~2.1.24]{MR2849819})

A matroid $M$ of rank $d\ge 2$ on the set $[n]$ is paving if and only if the collection of hyperplanes $\mathcal{H}$ of $M$, considered as a hypergraph with the vertex set $[n]$, satisfies the following conditions:

\begin{itemize}
    \item Every edge of $\mathcal{H}$ has cardinality at least $d-1$.

    \item Every subset of $[n]$ of cardinality $d-1$ is contained in exactly one edge of $\mathcal{H}$.
\end{itemize}

\end{definition}

It is worth mentioning that this definition has recently been extended to introduce a new class of matroids known as elementary split matroids by \citet[Theorem~6]{MR4499860}. The collection of hyperplanes of cardinality $d-1$ of a paving matroid is the collection of all subsets of $[n]$ of cardinality $d-1$ that is not contained in any hyperplane of cardinality greater than $d-1$. Therefore, the collection of hyperplanes of a paving matroid is determined by the collection of all hyperplanes of cardinality greater than $d-1$, or equivalently dependent hyperplanes. The following definition describes paving matroids in terms of dependent hyperplanes.

\begin{definition}\label{def}

A matroid $M$ of rank $d\ge 2$ on the set $[n]$ is paving if and only if the collection of dependent hyperplanes $\mathcal{H}$ of $M$ is a hypergraph with the set of vertices $[n]$, satisfying the following conditions:

\begin{itemize}
    \item Every edge of $\mathcal{H}$ has cardinality at least $d$.

    \item For any two distinct edges $A,B$ of $\mathcal{H}$, the cardinality of the intersection $A\cap B$ is at most $d-2$.
\end{itemize}
    
\end{definition}

A paving matroid will be identified with the hypergraph formed by its dependent hyperplanes, which is always a clutter, together with its rank $d$. This clutter will just be called the clutter of the paving matroid. We use the terms ``edge'' and ``dependent hyperplane'' interchangeably for a paving matroid. Moreover, a clutter that is the clutter of some paving matroid of rank $d$ will be called a $d$-\emph{paving clutter}. Dependent hyperplanes are precisely the split flacets of a paving matroid. Thus, a paving matroid of rank $d$ is a positroid (resp.\ interval positroid) if and only if its clutter is a CA $d$-paving clutter (resp.\ an interval $d$-paving clutter). We will say an edge of the clutter of a paving matroid is \emph{small} if its cardinality is $d$, or equivalently if it is the circuit-hyperplane of the paving matroid. Any edge that is not small will be referred to as a \emph{large} edge.

It is clear that a partial clutter, obtained by removing a dependent hyperplane of a given paving matroid, is also the set of dependent hyperplanes of a new paving matroid. This operation is called the \emph{relaxation} of a \emph{stressed hyperplane} introduced by \citet{10.1093/imrn/rnac270} which generalizes the relaxation of a circuit-hyperplane. Since a partial clutter of an interval clutter remains an interval clutter, the class of paving interval positroids is closed under taking relaxations.

Notice that a paving matroid contains a loop if and only if the rank of the matroid is $1$. And a paving matroid contains a coloop $x$ only if it is either the uniform matroid $U_{n,n}$ or $[n]\setminus\{x\}$ is the only dependent hyperplane. These cases will not appear since every paving matroid we consider will be of rank at least $2$ with at least three edges. Finally, the following lemma explains the effects of deletion and contraction operations on the clutter of a paving matroid.

\begin{lemma}\label{lem:delcon}

Let $\mathcal{H}=([n],E)$ be the clutter of a paving matroid $M$ of rank $d\ge 2$. Suppose that $x$ is not a coloop. The clutter of the paving matroid $M-x$ is $$\{H \mid x\notin H, H\in E\}\cup\{H-x \mid x\in H, H\in E, |H|>d\}.$$ And the clutter of the contraction $M/x$ is $$\{H-x \mid x\in H, H\in E\}.$$

\begin{proof}

Bases of $M-x$ are bases of $M$ not containing $x$. Therefore, a hyperplane of the deletion $M-x$ is a subset $H\setminus\{x\}$ where $H$ is a hyperplane of $M$. Since edges of the clutter of $M$ are dependent hyperplanes, the clutter of $M-x$ is $\{H\mid x\notin H, H\in E\}\cup\{H-x\mid x\in H, H\in E, |H|>d\}$. A hyperplane of the contraction $M/x$ is a subset $H\setminus\{x\}$ where $H$ is a hyperplane of $M$ containing $x$. Therefore, the clutter of $M/x$ is $\{H-x\mid x\in H, H\in E\}$.
    
\end{proof}
    
\end{lemma}

Therefore, a $d$-paving clutter $(V,E)$ is an excluded minor of interval positroids if and only if it is not an interval clutter and, for every $x\in V$, both $$\{H\mid x\notin H, H\in E\}\cup\{H-x\mid x\in H, H\in E, |H|>d\}$$ and $$\{H-x\mid x\in H, H\in E\}$$ are interval clutters. Moreover, if $M$ is an excluded minor of interval positroids and some relaxation $M'$ of $M$ is not an interval positroid, then $M'$ is also an excluded minor of interval positroids.

\section{Relaxation-minimal paving excluded minors of interval positroids}\label{sec:rel}

\begin{figure}
  \centering
\subcaptionbox{$O_{a_1,\ldots,a_n;d}$}{
    \begin{tikzpicture}
      \draw[thick] (30:2) ellipse [x radius=1.5, y radius=0.4, rotate=-60];
      \draw[thick] (90:2) ellipse [x radius=1.5, y radius=0.4];
      \draw[thick] (150:2) ellipse [x radius=1.5, y radius=0.4, rotate=60];
      \draw[thick] (210:2) ellipse [x radius=1.5, y radius=0.4, rotate=120];
      \draw[thick] (270:2) ellipse [x radius=1.5, y radius=0.4];
      \draw (330:2) node[circle, fill=black, inner sep=1.5pt] {};
      \draw (345:2) node[circle, fill=black, inner sep=1.5pt] {};
      \draw (315:2) node[circle, fill=black, inner sep=1.5pt] {};
      \node[] at (60:2.3) {$a_1$};
      \node[] at (120:2.3) {$a_2$};
      \node[] at (180:2.3) {$a_3$};
      \node[] at (240:2.3) {$a_4$};
    \end{tikzpicture}
  }
\subcaptionbox{$\Sh_{a,b,c;d}$}{
    \begin{tikzpicture}
    \draw[thick] (270:2) ellipse [x radius=3, y radius=0.5] {};
    \draw[thick] (180:2) ellipse [x radius=0.5, y radius=2.4] {};
    \draw[thick] (0,0) ellipse [x radius=0.5, y radius=2.5] {};
    \draw[thick] (0:2) ellipse [x radius=0.5, y radius=2.4] {};
    \node[] at (270:3) {D};
    \node[] at (180:3) {A};
    \node[] at (90:3) {B};
    \node[] at (0:3) {C};
    \node[] at (-2,-2) {$a$};
    \node[] at (0,-2) {$b$};
    \node[] at (2,-2) {$c$};
    \end{tikzpicture}
  }
  \subcaptionbox{$Y_{a,b;d}$}{
  \begin{tikzpicture}
    \draw[thick] (30:1) circle [radius=1.5]; {};
    \draw[thick] (150:1) circle [radius=1.5]; {};
    \draw[thick] (270:1) circle [radius=1.5]; {};
    \draw[very thick] (270:2) node[circle, fill=black, inner sep=1.5pt] {};
    \draw[very thick] (0,0) node[circle, fill=black, inner sep=1.5pt] {};
    \node[] at (0,1) {$\emptyset$};
    \node[] at (210:1) {$a$};
    \node[] at (30:1.7) {$d-b-1$};
    \node[] at (330:1) {$b$};
    \node[] at (150:3) {$A$};
    \node[] at (150:1.7) {$d-a-1$};
    \node[] at (270:3) {$C$};
    \node[] at (30:3) {$B$};
  \end{tikzpicture}
  }
\caption{Relaxation-minimal list of paving excluded minors of positroids}
\label{fig:relmin}
\end{figure}
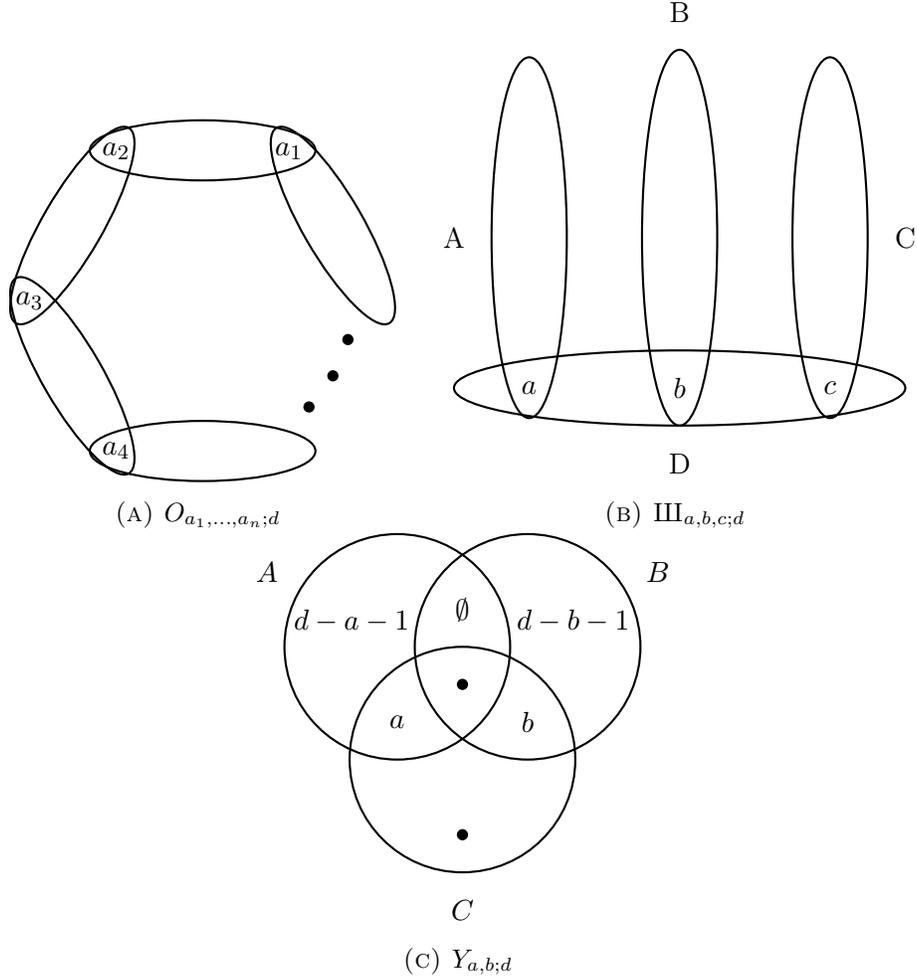

In this section, we show that every paving excluded minor of interval positroids can be reduced to relatively simple paving excluded minor of interval positroids through the relaxation of dependent hyperplanes. These ``relaxation-minimal'' paving excluded minors of interval positroids are illustrated in Figure \ref{fig:relmin}.

We first state the following fact.

\begin{proposition}

Any paving matroid of rank $2$ is an interval positroid.

\begin{proof}

The clutter of a paving matroid of rank $2$ is a collection of pairwise disjoint edges since any intersection of a pair of distinct edges has maximal size $d-2=0$. The result follows since a pairwise disjoint collection of edges forms an interval clutter.

\end{proof}
    
\end{proposition}

Consequently, every paving excluded minor of interval positroids must have a rank of at least $3$. The following theorem provides a characterization of relaxation-minimal paving excluded minors of interval positroids. 

\begin{lemma}\label{lem:noiso}

The clutter of any paving excluded minor of interval positroids does not contain an isolated point and any large edge contains at most one element exclusively.

\begin{proof}

Let $M$ be a paving matroid that is also an excluded minor of interval positroids. Then, the clutter of $M$ is not an interval clutter and proper deletion of $M$ is a paving matroid whose clutter is an interval clutter. Therefore, it is clear that the clutter of $M$ does not have an isolated point. Suppose that a large edge $A$ of $M$ contains more than one element exclusively. Choose $x\in A$ that is not contained in any other edge of $M$. Then, the collection of edges of the deletion $M-x$ consists of all edges of $M$ not containing $x$ and $A\setminus\{x\}$. This implies that $M-x$ is not an interval positroid. So, a large edge cannot contains more than one element exclusively.
    
\end{proof}
    
\end{lemma}

\begin{theorem}\label{thm:relmin}

Let $M$ be a paving excluded minor of interval positroids of rank $d\ge 3$. Let $\mathcal{H}=(V,E)$ be the clutter of $M$. Then, some relaxation of $M$, or equivalently, some partial clutter of $\mathcal{H}$ falls into one of the three types of non-interval clutters described below:

\begin{enumerate}
    \item For $m\ge 3$, let $1\le a_1,\ldots,a_m\le d-2$ be positive integers such that there are no three consecutive indices $i,i+1,i+2 \mod m$ satisfying two inequalities $$a_i+a_{i+1}>d,\ a_{i+1}+a_{i+2}>d.$$ The $d$-paving clutter $O_{a_1,\ldots,a_m;d}$ is defined by the clutter with $m$ edges $A_1,\ldots,A_m$ satisfying the following conditions:
        \begin{enumerate}
            \item No vertex is isolated.
        
            \item Every vertex is contained in at most two edges.

            \item $|A_i\cap A_{i+1}|=a_i$ where indices are taken modulo $m$.

            \item Each vertex in an edge with a cardinality greater than $d$ is also included in another edge with a cardinality of $d$.
        \end{enumerate}

    \item For nonnegative integers $0\le a,b\le d-3$, $Y_{a,b;d}$ is the $d$-paving clutter with three edges $A,B,C$ satisfying the following conditions:
        \begin{enumerate}
            \item No vertex is isolated.
        
            \item 
            \begin{align*}
                |A| = |B| &= d, \quad |C| \geq d \\
    |A \cap C \setminus B| &= a, \quad |C \cap B \setminus A| = b, \quad |A \cap B \setminus C| = 0 \\
    |C \setminus (A \cup B)| &> 0, \quad |A \cap B \cap C| = 1
            \end{align*}
            
            \item $|C|>d$ only if $C\setminus(A\cup B)$ contains exactly one element.
        \end{enumerate}

    \item For positive integers $1\le a,b,c\le d-2$, $\Sh_{a,b,c;d}$ is the $d$-paving clutter with four edges $A,B,C,D$ satisfying the following conditions:
        \begin{enumerate}
            \item No vertex is isolated.

            \item \begin{align*}
                |A|&=|B|=|C|=d, \quad |D|\ge d \\
                |A\cap D|&=a, \quad |B\cap D|=b, \quad |C\cap D|=c.
            \end{align*}
        
            \item $A,B,C$ are pairwise disjoint.

            \item $|D|>d$ only if $a+b+c=|D|$.
        \end{enumerate}
\end{enumerate}

\begin{proof}

Every non-interval partial clutter of $\mathcal{H}$ is also an excluded minor of interval positroids. Thus, a non-interval partial clutter cannot have an isolated vertex by Lemma \ref{lem:noiso}.

Since $\mathcal{H}$ is not interval, there is a partial clutter of $\mathcal{H}$ that is one of five non-interval clutter listed in Proposition \ref{prop:int}. Suppose that the clutter of $M$ has a partial clutter of type 1. So, there are three edges $A,B,C$ such that $(A\cap B)\setminus C, (B\cap C)\setminus A$ and $(A\cap C)\setminus B$ are all nonempty. The contraction of the set $A\cap B\cap C$ is a paving matroid with a partial clutter consisting of three edges $A\setminus (A\cap B\cap C), B\setminus (A\cap B\cap C)$ and $C\setminus (A\cap B\cap C)$ which remains a non-interval clutter of type 1. Therefore, $M$ is not an excluded minor of interval positroids unless $A\cap B\cap C$ is empty. If the edge $A$ is large and $A\setminus (B\cup C)$ is nonempty, then for $x\in A\setminus (B\cup C)$, the deletion $M-x$ has three edges $A\setminus\{x\}, B,C$ which is again a non-interval clutter of type 1. Thus, $A\subset B\cup C$ if $A$ is large. If both $A,B$ are large, then $C$ contains $(A\cup B)\setminus(A\cap B)$ by the previous argument. Then, either $|A\cap B|>1$ or $2d\le |A\cup B|-|A\cap B|=|C|$. If the latter case is true, then $A,B,C$ are all large. And some intersections of two distinct edges must contain more than one element. So, for both cases, two large edges share more than one element. Let $x$ be an element in the intersection of two large edges sharing more than one element. The deletion $M-x$ contains the partial clutter with three edges $A\setminus\{x\},B\setminus\{x\},C\setminus\{x\}$ which remains a non-interval clutter of type 1. Therefore, at most one edge can be large, and the partial clutter $(V,\{A,B,C\})$ is the clutter $O_{a_1,a_2,a_3;d}$ with $a_1=|A\cap B|, a_2=|B\cap C|, a_3=|C\cap A|$.

Let $A,B,C\in E$ form a partial clutter of type 2. So, $A, B, C$ share a common point, and each of these contains an element exclusively. If all three sets $(A\cap B)\setminus C$, $(B\cap C)\setminus A$ and $(A\cap C)\setminus B$ are nonempty, then the contraction $M/(A\cap B\cap C)$ contains a partial clutter with three edges $A\setminus (A\cap B\cap C)$, $B\setminus (A\cap B\cap C)$, and $C\setminus (A\cap B\cap C)$ which is a non-interval clutter of type $1$. By symmetry, we may assume that $(A\cap B)\setminus C$ is empty. If there is more than one element in the intersection $A\cap B\cap C$, then the contraction of an element $x$ contained in $A\cap B\cap C$ will be a paving matroid with the partial clutter $A\setminus\{x\}$, $B\setminus\{x\}$, $C\setminus\{x\}$ which remains a non-interval clutter of type $2$. Hence, $|A\cap B\cap C|=1$. Since the partial clutter with three edges $A,B,C$ is a $d$-paving clutter, the following inequality holds: $$|A\setminus (B\cup C)|=|A|-|A\cap B|\ge 2.$$ This implies that $A$ must be a small edge by Lemma \ref{lem:noiso}. By symmetry, the edge $B$ is also small. And $|C|>d$ only if $|B\setminus(A\cup C)|=1$ also by Lemma \ref{lem:noiso}. Thus, the partial clutter $(V,\{A,B,C\})$ is $Y_{a,b;d}$ where $a=|(A\cap C)\setminus B|$, $b=|(C\cap B)\setminus A|$.

Suppose that there are four edges $A,B,C,D$ of $\mathcal{H}$ that form a non-interval clutter of type 3. Contracting common points of $A,B,C,D$, the partial clutter remains non-interval since it contains the non-interval CA subhypergraph $\{12,23,34,41\}$.

If there is a non-interval partial clutter of type $4$, the condition $$1\leq |A\cap D|,\ |B\cap D|,\ |C\cap D|\leq d-2$$ immediately follows from the fact that $A, B, C, D$ form a $d$-paving clutter. The inequality $|A\setminus (B\cup C\cup D)|=|A\setminus D|>1$ implies that $A$ is small by Lemma \ref{lem:noiso}. By symmetry, the other two edges $B,C$ are also small. If $D$ is large and $D$ contains some element $x$ not contained in $A\cup B\cup C$, then the deletion $M-x$ will contain the partial clutter $(V,\{A,B,C,D\setminus\{x\}\})$ which remains the non-interval clutter of type 4. Therefore, $D$ can be possibly large only if every element of $D$ is also contained in one of $A,B,C$. Hence, the partial clutter $A,B,C,D$ is the clutter $\Sh_{a,b,c;d}$.

Finally, for some $m\ge 4$, suppose that there are $A_1,\ldots, A_m\in E$ such that $A_i\cap A_j\neq\emptyset$ if and only if $i=j\pm 1 \mod m$ which is a non-interval clutter of type 5. Let $x$ be a vertex that is contained in only large edges. Since every vertex is contained in at most two edges, it is either contained in a large edge exclusively or it is contained in the intersection of two large edges. If $x$ is contained in only one large edge $A_i$, the deletion $M-x$ contains the non-interval partial clutter $A_1,\ldots,A_i\setminus\{x\},\ldots,A_m$ of type 5. This is a contradiction since $M$ is an excluded minor of interval positroids. Now, suppose that $x$ is contained in the intersection $A_i\cap A_{i+1}$ and $A_i, A_{i+1}$ are both large. If $|A_i\cap A_{i+1}|>1$, then the deletion $M-x$ contains the partial clutter $A_1,\ldots, A_i\setminus\{x\}, A_{i+1}\setminus\{x\},\ldots,A_m$ which is a non-interval clutter of type 5. So, the intersection $A_i\cap A_{i+1}$ contains exactly one element. Since it is a $d$-paving clutter, the inequality $|A_i\cap A_{i-1}|\le d-2$ holds, which implies that the set $A_i\setminus(A_{i-1}\cup A_{i+1})$ is nonempty. This goes back to the previous case and yields a contradiction. Therefore, any vertex contained in a large edge is contained in another small edge. So the partial clutter $(V,\{A_1,\ldots,A_m\})$ is $O_{a_1,\ldots,a_m;d}$ for $a_i=|A_i\cap A_{i+1}|$.

\end{proof}

\end{theorem}

\begin{remark}

Every paving excluded minor of interval positroids characterized in \ref{thm:relmin} is ``relaxation-minimal'', meaning that further relaxation leads to an interval positroid.
    
\end{remark}

\section{Paving excluded minors of interval positroids that are positroids}\label{sec:pos}

Suppose that $M$ is a paving positroid that is an excluded minor of interval positroids. Then, the clutter of $M$ is a CA clutter. We provide a criterion for testing whether a given paving positroid $M$ is an excluded minor of interval positroids.

\begin{definition}

A non-interval CA clutter with three edges is a triangle.
    
\end{definition}

\begin{definition}

A clutter with three edges $A,B,C$ will be called a triangle with a common point if $A\cap B\cap C\neq\emptyset$, $A\cap C\not\subset B$, and $B\subsetneq A\cup C$.
    
\end{definition}

It is clear that a clutter $\mathcal{H}$ with three edges $A,B,C$ with a common point $x\in A\cap B\cap C$ is CA if and only if at least one of $A,B,C$ contains no element exclusively. If $B$ is contained in $A\cup C$, then $\mathcal{H}$ is non-interval if and only if $A\cap C\not\subset B$. Hence, a clutter with three edges $A, B, C$, and the nonempty intersection $A \cap B \cap C$ is a triangle if and only if it is a triangle with a common point, explaining the terminology.

\begin{proposition}\label{prop:posex}

A paving positroid $M$ of rank $d$ is an excluded minor of interval positroids if and only if the clutter $\mathcal{H}=(V,E)$ of $M$ satisfies the following conditions:

\begin{enumerate}
    \item $\mathcal{H}$ is a non-interval, CA $d$-paving clutter.

    \item If a partial clutter of $\mathcal{H}$ with three edges is a non-interval CA clutter, then it is not a triangle with a common point and at most one of the three edges is large.

    \item Let $A,B$ be two (not necessarily distinct) large edges with the nonempty intersection $A\cap B$. If $C, D$ are edges overlapping with $A\cup B$ such that $(A\cup B)\cap C, (A\cup B)\cap D$ are incomparable, then $A\cup B\subset C\cup D$.       
\end{enumerate}

\begin{proof}

The first condition is clear since $\mathcal{H}$ is the clutter of a paving positroid that is not an interval positroid. Since a triangle with two large edges or a common point is not an excluded minor of interval positroids by Theorem \ref{thm:relmin}, it is necessary that every partial clutter of $\mathcal{H}$ that is a triangle has at most one large edge without a common point. Now, suppose that $A,B$ are large edges with a nonempty intersection. Fix an arc-ordering $f:V\to [n]$ and we identify each vertex and edge with its image under $f$. Then, the union $A\cup B$ is a cyclic interval. Let $C,D$ are edges overlapping with the cyclic interval $A\cup B$. If two intersections $(A\cup B)\cap C, (A\cup B)\cap D$ are incomparable, then both endpoints of $A\cup B$ is contained in $C\cup D$. Assume, for a contradiction, that $C\cup D$ does not contain $A\cup B$. Then, the following two cases are possible.

\begin{enumerate}
    \item There is an element $x\in A$ not contained in $B\cup C\cup D$ (or $x\in B$ not contained in $A\cup C\cup D$ by symmetry).

    \item There is an element $x\in A\cup B$ not contained in $C\cup D$ and $|A\cup B|\ge 2$.
\end{enumerate}

Suppose that the first condition does not hold. If $|A\cap B|=1$, then $|A\setminus B|\ge d, |B\setminus A|\ge d$ since $A,B$ are large. The clutter $\mathcal{H}$ is a paving $d$-clutter and therefore $|(A\cup B)\cap C|,\ |(A\cup B)\cap D|\le d-2$ which implies that $(A\setminus B)\cup (B\setminus A)$ can not be contained in $C\cup D$. This is a contradiction since the first condition does not hold. Thus, either the first or the second condition holds. Now, $\mathcal{H}$ remains non-interval after deleting an element $x\in A$ not contained in $B\cup C\cup D$ or $x\in A\cup B$ not contained in $C\cup D$ when $|A\cap B|\ge 2$.

Conversely, suppose that the contraction $M/x$ is a non-interval, CA hypergraph for some $x\in V$. By Lemma \ref{lem:delcon}, the edges of the clutter of $M/x$ is the collection of all sets $H\setminus\{x\}$ where $H\in E$ and $x\in H$. Therefore, the collection $\mathcal{H}_x$ of all edges of $\mathcal{H}$ containing $x$ is a non-interval and CA partial clutter. This implies that there are more than two edges in the clutter $\mathcal{H}_x$. By Corollary \ref{cor:ca}, the clutter $\mathcal{H}_x$ contains a non-interval partial clutter of type 1 or 5. But edges of the type $5$ clutters can not share a vertex in common. Hence, there is a partial clutter of type 1 with three edges in $\mathcal{H}_x$ that is a CA clutter sharing a point $x$. This is a triangle with a common point.

Now, assume that the clutter $\mathcal{H}-x$ of the deletion $M-x$ is not an interval clutter for some $x\in [n]$. Again, by Corollary \ref{cor:ca}, this is equivalent to some partial clutter of $\mathcal{H}-x$ is a non-interval clutter of type 1 or 5. Then, there are the following four possibilities:

\begin{enumerate}
    \item $x$ is contained in some three large edges covering the vertex set $V$ and forming a non-interval clutter of type 1.

    \item There is a non-interval partial clutter of type 1 such that $x$ is exactly contained in two adjacent large edges of the partial clutter.

    \item There is a non-interval partial clutter of type 5 such that $x$ is exactly contained in two adjacent large edges of the partial clutter.

    \item There is a non-interval partial clutter of type 1 or 5 such that $x$ is exactly contained in one large edge of the partial clutter.
\end{enumerate}

The first case is impossible if there is no partial clutter which is a triangle with a common point. The second case is impossible if there is no partial clutter which is a triangle with two large edges. The third and fourth cases cannot occur unless the large edges $A$ and $B$ containing $x$ violate the third condition. Therefore, the three given conditions are sufficient.
    
\end{proof}
    
\end{proposition}

\section{Paving excluded minors of interval positroids with four dependent hyperplanes}\label{sec:four}

In this section, we classify all paving, non-positroid excluded minors of interval positroids with four dependent hyperplanes. By Theorem \ref{thm:relmin}, it is either the clutter $\Sh_{a,b,c;d}$, or a $d$-paving clutter obtained by adding one edge to $Y_{a,b;d}$ or $O_{a_1,a_2,a_3;d}$.

\begin{theorem}\label{thm:Y4}

Let $M$ be a paving matroid of rank $d$ with exactly four edges $A,B,C_1,C_2$. If the partial clutter $(V,E=\{A,B,C_1\})$ is $Y_{a_1, b_1; d}$ with the common point $v$, then $M$ is an excluded minor of interval positroids if and only if one of the following statements hold:

\begin{enumerate}
    \item $|C_1|=d+1$ and $C_2=(A\setminus C_1)\cup (C_1\setminus(A\cup B))\cup (B\setminus C_1)$.

    \item The partial clutter $(V,E=\{A,B,C_2\})$ is $Y_{a_2, b_2; d}$ where $|(A\cap C_2)\setminus\{v\}|=a_2, |(B\cap C_2)\setminus\{v\}|=b_2$ with the same common point $v$, and the clutter $(V\setminus\{v\},\{A\setminus\{v\},B\setminus\{v\},C_1\setminus\{v\},C_2\setminus\{v\}\})$ is an interval $(d-1)$-paving clutter.
\end{enumerate}

These clutters will be denoted by $Y_{a_1,b_1;d}', Y_{a_1,a_2,b_1,b_2;d}$ respectively.

\begin{proof}
Checking whether these clutters are indeed excluded minors is omitted.

~
\paragraph{\textbf{Case 1:} $v\notin C_2$.}$ $\newline

If $C_2$ is contained in $A\cup C_1$, then the partial clutter with three edges $A,C_1,C_2$ is a triangle with a common point. Hence, the intersection $C_2\cap(A\setminus C_1)$ is nonempty. By symmetry, the intersection $C_2\cap (B\setminus C_1)$ is also nonempty. If there is an element $x\in C_1\setminus(A\cup B)$ not contained in the edge $C_2$, then the matroid $M-x$ has three edges $A,B,C_2$ which is a non-interval clutter of type 1. Therefore, the set $C_1\setminus(A\cup B)$ must be contained in $C_2$. For the same reason, both $A\setminus C_1$ and $B\setminus C_1$ are contained in $C_2$. Now, if there is an element $x$ contained in $A\cap C_1\cap C_2$, then the contraction $M/x$ is a non-interval clutter of type 1. Therefore, $A\cap C_1\cap C_2$ is empty. By symmetry, $B\cap C_1\cap C_2$ is also empty. Hence, $C_2$ must be $(A\setminus C_1)\cup (C_1\setminus(A\cup B))\cup (B\setminus C_1)$.

If $|C_2|>d$, for $x\in A\setminus C_1$, $M/x$ has three edges $B,C_1,C_2\setminus\{x\}$ which is a non-interval clutter of type 1. Hence, $|C_2|=d$. If $|A|=|B|=|C_1|=d$, then $a+b=|C_1\cap(A\cup B)|\le d-2$ since $C_1\setminus (A\cup B)$ is nonempty. But

\begin{align*}
    |C_2| &= d \\
        &= |A\setminus C_1| + |C_1\setminus(A\cup B)| + |B\setminus C_1| \\
        &= (d-a-1) + (d-a-b-1) + (d-b-1) \\
        &= 3(d-1) - 2a - 2b
\end{align*} implies that $2d=2a+2b+3$ which is a contradiction since $a+b\le d-2$. Finally, assume that the edge $C_1$ is large. This implies that $|C_1\setminus(A\cup B)|=1$ and $a+b=|C_1|-2\ge d-1$ by Theorem \ref{thm:relmin}. But, we have

\begin{align*}
    |C_2| &= d \\
        &= |A\setminus C_1| + |C_1\setminus(A\cup B)| + |B\setminus C_1| \\
        &= (d-a-1) + 1 + (d-b-1) \\
        &= 2d - a - b-1.
\end{align*} This implies that $a+b=d-1$ and $|C_1|=d+1$. This gives us the first case in the list.

~
\paragraph{\textbf{Case 2:} $v\in C_2$.}$ $\newline

Let $X'$ be the set $X\setminus\{v\}$ for edges $X=A,B,C_1,C_2$. Since $v\in C_2$, the contraction $M/v$ has four edges $A',B',C_1',C_2'$ which must be an interval clutter. If both $A'\setminus C_1'$ and $B'\setminus C_1'$ share an element with $C_2'$, then $C_2'$ must contain $C_1'$ which is impossible since $(V,\{A,B,C_1,C_2\})$ is a clutter. Hence, $C_2'$ is properly contained in either $A'\cup C_1'$ or $B'\cup C_1'$. By symmetry, we may assume that $C_2'\subsetneq B'\cup C_1'$. To avoid a triangle with a common point as a partial clutter, the intersection $B'\cap C_1'$ is contained in $C_2'$. If there is an element $x\in C_1\setminus(A\cup B)$ not contained in $C_2$, then the clutter of the deletion $M-(C_1\setminus(A\cup B\cup C_2))$ is the clutter $Y_{0,|C_2'\cap B'|;d}$ with three edges $A,B,C_2$ which is not interval. Hence, $C_2'$ contains the set $C_1\setminus(A\cup B)$. Thus, the partial clutter $A,B,C_2$ is the clutter $Y_{|A'\cap C_2'|,|C_2'\cap B'|;d}$ with $C_1\setminus (A\cup B)=C_2\setminus(A\cup B)$. 
    
\end{proof}
    
\end{theorem}

Before classifying paving excluded minors of interval positroids obtained by adding an edge to the clutter $O_{a_1,a_2,a_3;d}$, we introduce the following notation for convenience. Edges of $O_{a_1,a_2,a_3;d}$ will be denoted by $A,B,C$ and let \begin{align*}
    S&=A\setminus(B\cup C), S'=B\cap C, \\
    T&=B\setminus(C\cup A), T'=C\cap A, \\
    U&=C\setminus(A\cup B), U'=A\cap B, \\
    a&=|S'|, b=|T'|, c=|U'|.
\end{align*}

Notice that if $|A|=|B|=|C|=d$, then $|S|=a+k$, $|T|=b+k$, and $|U|=c+k$ for some $k\ge -\min\{a,b,c\}$. Throughout the paper, the variable $k$ will consistently represent this difference in cardinality.

\begin{theorem}\label{thm:O4}

Let $M$ be a paving excluded minor of interval positroids of rank $d$ with exactly four dependent hyperplanes $A,B,C,D$. Suppose that the partial clutter $(V,E=\{A,B,C\})$ is $O_{a,b,c;d}$. If $M$ is not a positroid, then $|A|=|B|=|C|=d$, and the fourth edge $D$ of the clutter $\mathcal{H}$ of $M$ is one of the following two possibilities:

\begin{enumerate}
    \item $k=0$ and $D=S\cup T\cup U$.

    \item $D=S'\cup T'\cup U'$ for $a,b,c,k$ satisfying $$2\le |S|=a+k\le a,\ 2\le |T|=b+k\le b,\ 2\le |U|=c+k\le c.$$    
\end{enumerate}

The first clutter will be denoted by $\Delta_{a,b,c}^1$ and the second clutter will be denoted by $\Delta_{a,b,c;k}^2$.

\begin{proof}

Verifying that these two clutters are indeed excluded minors of interval positroids is straightforward. So, we focus on proving that these are only possibilities.

~
\paragraph{\textbf{Case 1:} No element of $\mathcal{H}$ is contained in three edges.}$ $\newline

Hence, $D$ is contained in the set $S\cup T\cup U$. If $A$ is large, then $B$ and $C$ are small and $A\subset B\cup C$ by Theorem \ref{thm:relmin}. Therefore, the cardinality of the set $S\cup T\cup U$ is $|B|+|C|-|A|-|S|<2d-|A|\le d-1$. This is a contradiction since $d\le |D|$. Therefore $A,B,C$ are small. Since $D\subseteq S\cup T\cup U$ and $d=a+b+c+k,\ |S\cup T\cup U|=a+b+c+3k$, $k$ must be some nonnegative integer. This implies that $S$, $T$, and $U$ are all nonempty. Since $\mathcal{H}$ is a clutter, $D$ is not contained in any of $S$, $T$, and $U$. By symmetry, we may assume that $S\cap D$ and $T\cap D$ are nonempty. If there is an element $x\in U$ not contained in $D$, then the clutter of the deletion $M-x$ has three edges $A$, $B$, and $D$ which is a non-interval clutter of type 1 since every pairwise intersection is nonempty. Therefore, it follows that $U \subset D$, and similarly, both $S$ and $T$ are also contained in $D$. Thus, we have shown that $D=S\cup T\cup U$. It remains to show that $D$ is small. Note that $D$ is large if and only if $k>0$. If $D$ is large, then for any element $x\in U$, the clutter of the deletion $M-x$ contain three edges $A$, $B$, and $D\setminus\{x\}$. The clutter with three edges $A$, $B$, and $D\setminus\{x\}$ is a non-interval clutter of type 1. Therefore, $M-x$ is not an interval positroid if $D$ is large. So, $k=0$ and $|S|=a, |T|=b, |U|=c$, and $D=S\cup T\cup U$.

~
\paragraph{\textbf{Case 2:} Some element $v\in S'$ is contained in $D$, and $D\cap T'$ and $D\cap U'$ are empty.}$ $\newline

~
\paragraph{Step 1: Three edges $A$, $B$, and $C$ are small.}$ $\newline

If $A$ is large, then $S=\emptyset$ which implies that $D$ is contained in $S'\cup T\cup U$. Since $B,C$ are small, $|S'|+|T|+|U|<d$. This is a contradiction because $D$ is contained in $S\cup S'\cup T\cup U$. If $B$ is large, then $B=S'\sqcup U'$. This implies that $1<|S'|, |U'|$ since $|S'|, |U'|\le d-2$. Thus, $D\cap S$ is nonempty since $D$ is not contained in $C$. For an element $x\in U'$, the clutter of the deletion $M-x$ consists of three edges $B\setminus\{x\}$, $C$, and $D$ sharing the element $v$. Elements of $U'\setminus\{x\}$, $T'$, and $D\cap S$ are exclusively contained in the edges $B\setminus\{x\}$, $C$, and $D$ respectively. Therefore, the clutter of $M-x$ is not an interval clutter. By symmetry, the edge $C$ can not be large either.

~
\paragraph{Step 2: $D=S\cup T\cup U\cup S'$.}$ $\newline

The intersection $D \cap S$ must be nonempty because if not, $D$ is contained in the set $B \cup C$, which implies that either $\mathcal{H}$ is interval, or the partial clutter with edges $B$, $C$, and $D$ forms a triangle with a common point. Then, the partial clutter obtained by removing the edge $B$ is non-interval. Therefore, $T$ must be entirely contained in $D$; otherwise, we could delete an element from $T\setminus D$ to obtain a proper minor with three edges $A, C, D$. By symmetry, $U$ is also contained in $D$. If there is an element $x\in S'$ not contained in $D$, then the clutter of the contraction $M/v$ consists of three edges $B\setminus\{v\}$, $C\setminus\{v\}$, and $D\setminus\{v\}$ with $((B\setminus\{v\})\cap (C\setminus\{v\}))\setminus (D\setminus\{v\})$ nonempty. This leads to the false conclusion that $M/v$ has three edges $B\setminus\{v\}$, $C\setminus\{v\}$, and $D\setminus\{v\}$ which form a non-interval clutter of type 1. So, $S'$ is also contained in $D$. Finally, if there is an element $x\in S$ not contained in $D$, then the clutter of the deletion $M-x$ has three edges $B,C,D$ which form a non-interval clutter of type 2 in Proposition \ref{prop:int}. Therefore, $D=S\cup T\cup U\cup S'$. The cardinality of $D$ is $|S|+|T|+|U|+|S'|=a+b+c+3k+1$. Since $d\le |D|$ and $d=a+b+c+k$, $k$ must be some nonnegative integer. Therefore, $D$ is large. Then, for any $x\in T$, the deletion $M-x$ has three edges $A,C,D\setminus\{x\}$ which is non-interval. Therefore, this case is impossible.

~
\paragraph{\textbf{Case 3:} $S'\cap D$ and $T'\cap D$ are nonempty and $U'\cap D$ is empty.}$ $\newline

Since the partial clutter $A$, $B$, and $D$ is a non-interval clutter of type 1, $U$ is contained in $D$. If $D$ is contained in $A\cup C$, then either $\mathcal{H}$ is an interval clutter or the partial clutter $A,C,D$ is a triangle with a common point. Hence, the intersection $D\cap T$ is nonempty. By symmetry, the intersection $D\cap S$ is also nonempty. Since the edge $C$ is not contained in $D$, either $S'\not\subset D$ or $T'\not\subset D$. If $S'\not\subset D$, then the contraction $M/v$ of an element $v\in D\cap S'$ has three edges $B\setminus\{v\}, C\setminus\{v\}, D\setminus\{v\}$ which form a non-interval clutter of type 1. Therefore, this case is impossible.

~
\paragraph{\textbf{Case 4:} $S'\cap D$, $T'\cap D$, and $U'\cap D$ are all nonempty.}$ $\newline

If some element of $S'$ is not contained in $D$, then the clutter of the contraction $M/x$ of an element $x\in S'\cap D$ is a non-interval clutter of type 1. Hence $S'$ is contained in the edge $D$. By symmetry, the set $S'\cup T'\cup U'$ is contained in the edge $D$.

~
\paragraph{Step 1: Three edges $A,B,C$ are small.}$ $\newline

If $A$ is large, then $A=T'\sqcup U'$. Therefore $A$ is contained in $D$ by the above observation. By symmetry, $A,B,C$ are all small.

~
\paragraph{Step 2: $D=S'\cup T'\cup U'$.}$ $\newline

 Since edges $A,B,C$ are not contained in the edge $D$, there are elements $x\in S, y\in T, z\in U$ that are not contained in $D$. Now, if there is an element in $S$ contained in $D$, then the clutter of the deletion $M-x$ consists of three edges $B, C, D$ sharing vertices in $S'$, and each of which contains an element exclusively. Thus, three edges $B,C,D$ form a non-interval clutter of type 2. Hence, $D$ must be $S'\cup T'\cup U'$.
 
~
\paragraph{Step 3: $2\le |S|\le a$, $2\le |T|\le b$, $2\le |U|\le c$.}$ $\newline
 
$|S|=a+k, |T|=b+k, |U|=c+k$ because $A,B,C$ are small. Since $d=a+b+c+k$ and $|D|=a+b+c$, $k$ is not positive. Moreover, the cardinality of the set $S=A\setminus D$ is at least $2$ since $\mathcal{H}$ is a $d$-paving clutter. By symmetry, the given inequalities hold.
    
\end{proof}
    
\end{theorem}

\begin{figure}
\captionsetup[subfigure]{labelformat=empty}
    \centering
    \subcaptionbox{$\Sh_{a,b,c;d}$}{
    \begin{tikzpicture}
    \draw[thick] (270:2) ellipse [x radius=3, y radius=0.5] {};
    \draw[thick] (180:2) ellipse [x radius=0.5, y radius=2.4] {};
    \draw[thick] (0,0) ellipse [x radius=0.5, y radius=2.5] {};
    \draw[thick] (0:2) ellipse [x radius=0.5, y radius=2.4] {};
    \node[] at (270:3) {D};
    \node[] at (180:3) {A};
    \node[] at (90:3) {B};
    \node[] at (0:3) {C};
    \node[] at (-2,-2) {$a$};
    \node[] at (0,-2) {$b$};
    \node[] at (2,-2) {$c$};
    \end{tikzpicture}
  }
  \subcaptionbox{$Y'_{a_1,b_1;d}$ for $a_1+b_1=d-1$}{
  \begin{tikzpicture}
    \draw[thick] (30:1) circle [radius=1.5]; {};
    \draw[thick] (150:1) circle [radius=1.5]; {};
    \draw[thick] (270:1) circle [radius=1.5]; {};
    \draw[very thick] (270:2) node[circle, fill=black, inner sep=1.5pt] {};
    \draw[very thick] (0,0) node[circle, fill=black, inner sep=1.5pt] {};
    \node[] at (0,1) {$\emptyset$};
    \node[] at (210:1) {$a_1$};
    \node[] at (150:1.7) {$d-a_1-1$};
    \node[] at (330:1) {$b_1$};
    \node[] at (150:3) {$A$};
    \node[] at (270:3) {$C$};
    \node[] at (30:3) {$B$};
    \node[] at (30:1.7) {$d-b_1-1$};
    \draw[ultra thick, draw=blue] (160:1.7) to[out=225, in=180, looseness=1.5] (265:2);
    \draw[ultra thick, draw=blue] (275:2) to[out=0, in=315, looseness=1.5] (20:1.7);
  \end{tikzpicture}
  }
  \subcaptionbox{$Y_{a_1,a_2,b_1,b_2;d}$}{
  \begin{tikzpicture}
    \draw[thick] (30:1) circle [radius=1.5]; {};
    \draw[blue, ultra thick] (30:1) ++(-65:1.5) arc [start angle=-65, end angle=0, radius=1.5];
    \draw[blue, ultra thick] (13:2.4) -- (30:0.7);
    \draw[blue, ultra thick] (30:0.7) arc [start angle=65, end angle=115, radius=1.5];
    \draw[blue, ultra thick] (150:0.7) arc [start angle=185, end angle=210, radius=1.5];
    \draw[blue, ultra thick] (210:1.7) -- (210:0.5);
    \draw[blue, ultra thick] (210:1.7) arc [start angle=-185, end angle=5, radius=1.5];
    \draw[thick] (150:1) circle [radius=1.5]; {};
    \draw[thick] (270:1) circle [radius=1.5]; {};
    \draw[very thick] (0,0) node[circle, fill=black, inner sep=1.5pt] {};
    \node[] at (0,1) {$\emptyset$};
    \node[] at (150:3) {$A$};
    \node[] at (210:1) {$a_1$};
    \node[] at (330:1) {$b_1$};
    \node[] at (270:3) {$C_1$};
    \node[] at (270:1.5) {$\neq\emptyset$};
    \node[] at (270:0.3) {$v$};
    \node[] at (30:3) {$B$};
  \end{tikzpicture}
  }
    \subcaptionbox{$\Delta_{a,b,c}^1$}{
    \begin{tikzpicture}
  \begin{scope}[rotate=300]
    \draw[thick] (30:1) ellipse [x radius=1.9, y radius=0.5] node[] (b) {};
  \end{scope}
  \begin{scope}[rotate=60]
    \draw[thick] (150:1) ellipse [x radius=1.9, y radius=0.5] node[] (c) {};
  \end{scope}
  \draw[thick] (270:2) ellipse [x radius=1.9, y radius=0.5] node[] (a) {};
  \draw[ultra thick, draw=blue, rounded corners=8pt] (a) -- (b) -- (c) -- (a);

  \node at (a) [] {a};
  \node at (a) [below, yshift=-0.5cm] {A};
  \node at (b) [] {b};
  \node at (b) [xshift=0.5cm, yshift=0.5cm] {B};
  \node at (c) [] {c};
  \node at (c) [xshift=-0.5cm, yshift=0.5cm] {C};
  \node at (90:0.9) {$a$};
  \node at (232:2.4) {$b$};
  \node at (308:2.4) {$c$};
\end{tikzpicture}}
\subcaptionbox{$\Delta_{a,b,c;k}^2$ for $k\le 0$}{
\begin{tikzpicture}
  \begin{scope}[rotate=300]
    \draw[thick] (30:1) ellipse [x radius=1.9, y radius=0.5] node[] (b) {};
  \end{scope}
  \begin{scope}[rotate=60]
    \draw[thick] (150:1) ellipse [x radius=1.9, y radius=0.5] node[] (c) {};
  \end{scope}
  \draw[thick] (270:2) ellipse [x radius=1.9, y radius=0.5] node[] (a) {};
  \node at (a) [] {$a+k$};
  \node at (a) [below, yshift=-0.5cm] {$A$};
  \node at (b) [] {$b+k$};
  \node at (b) [xshift=0.5cm, yshift=0.5cm] {$B$};
  \node at (c) [] {$c+k$};
  \node at (c) [xshift=-0.5cm, yshift=0.5cm] {$C$};
  \node at (90:0.9) {a};
  \node at (232:2.4) {b};
  \node at (308:2.4) {c};
  \draw[ultra thick, draw=blue] (100:0.9) to[out=180, in=180, looseness=1.5] (230:2.5);
  \draw[ultra thick, draw=blue] (80:0.9) to[out=0, in=0, looseness=1.5] (310:2.5);
\end{tikzpicture}
}
    \caption{Non-positroid paving excluded minors of interval positroids with four edges}
    \label{fig:ex4}
\end{figure}
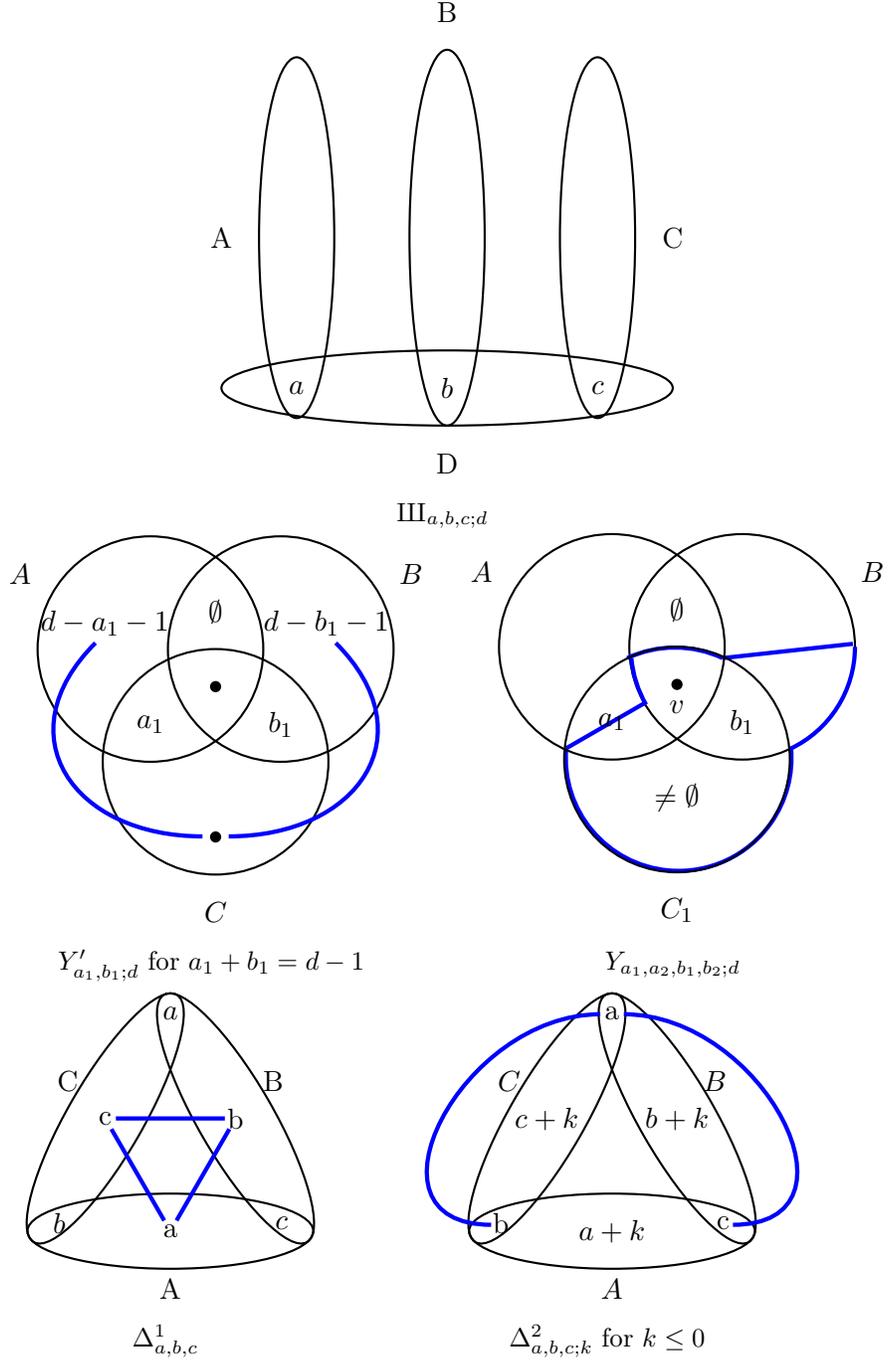

In summary, Figure \ref{fig:ex4} illustrates all non-positroid, paving matroids with four edges that are excluded minors of interval positroids. The fourth edge is indicated by blue arcs and segments.

\section{Paving, non-positroid excluded minors of interval positroids with more than five edges}\label{sec:five}

In this section, we classify all paving excluded minors of interval positroids with more than four edges. Again, we focus on paving excluded minors of interval positroids that are not positroids.

\subsection{Local compatibility implies global compatibility for CA clutters}

Let $\mathcal{H}=(V,E=\{A_1,\ldots,A_m\})$ be a non-interval CA clutter with every vertex is contained in at most two edges throughout this subsection. By Lemma \ref{lem:arcorder}, every arc ordering of $\mathcal{H}$ maps each intersection $A_i\cap A_{i+1}$ and each set of isolated points $A_i\setminus(A_{i-1}\cup A_{i+1})$ to disjoint cyclic intervals. These cyclic intervals will be called the ``blocks'' of $\mathcal{H}$. These blocks appear in the order $A_1\cap A_m, A_1\setminus(A_m\cup A_2), A_1\cap A_2,\ldots,A_m\cap A_{m-1}, A_m\setminus(A_1\cup A_{m-1})$ for every arc ordering. Let $T_1,T_2,\ldots,T_{2m}$ be these blocks in the given order. Let $\mathcal{H}+B$ denote the hypergraph obtained by adding a new edge $B$ to the hypergraph $\mathcal{H}$. We state the necessary and sufficient condition for $\mathcal{H}+B$ to be a CA clutter in the following lemma.

\begin{lemma}\label{lem:intuni}

Let $\mathcal{H}=(V,E=\{A_1,\ldots,A_m\})$ be a non-interval CA clutter such that every vertex is contained in at most two edges. The hypergraph $\mathcal{H}+B$, obtained by adjoining another edge $B$, is also a CA clutter if and only if

\begin{enumerate}
    \item for some unique $i$, $A_i\cap A_{i+1}\subset B\subset A_i\cup A_{i+1}$;

    \item the set of indices $I:=\{i\in [2m]\mid T_i\cap B\neq\emptyset\}$ is a cyclic interval of $[2m]$ and $B$ and $T_i$ overlap only if $i$ is one of the two endpoints of $I$.
\end{enumerate}

\begin{proof}

The sufficiency of two conditions is clear since every arc-ordering of $\mathcal{H}$ maps each block $T_1,\ldots,T_{2m}$ to disjoint cyclic intervals in the given order by Lemma \ref{lem:arcorder}. The necessity of the second condition also follows from this fact. Now, it remains to prove the necessity of the first condition.

Suppose that an arc ordering of $\mathcal{H}$ is given such that $B$ is a cyclic interval. If $B$ contains four consecutive blocks, then $B$ contains $A_i$ for some $i$ which is a contradiction. If the cyclic interval $B$ contains no block, then $B$ is contained in some $A_i$ which is also a contradiction.

Now, if the cyclic interval $B$ contains exactly one block $T_j$, then $B$ is contained in the union of three consecutive blocks $T_{j-1}\cup T_j\cup T_{j+1}$. Therefore $T_j$ must be some intersection $A_i\cap A_{i+1}$ because otherwise $B$ is contained in some edge $A_i$. If $B$ contains exactly two consecutive blocks, then one of the two blocks must be the intersection of two edges $A_i\cap A_{i+1}$. Finally, if $B$ contains three consecutive blocks, then the middle block must be some intersection $A_i\cap A_{i+1}$ since $B$ does not contain some $A_i$. For all these cases, it is clear that $B$ is contained in the union $A_i\cup A_{i+1}$.

\end{proof}
    
\end{lemma}

The next lemma characterizes three possible cases of adding two edges to construct a new CA clutter from $\mathcal{H}$.

\begin{lemma}\label{lem:incomp}

Let $\mathcal{H}_{12}:=(V,E\cup\{B_1,B_2\})$ be a clutter obtained by adding two edges to $\mathcal{H}$. Let $\mathcal{H}_i:=(V,E\cup\{B_i\})$ be the clutter obtained by only adding one edge $B_i$. Suppose that $\mathcal{H}_i$ is CA for $i=1,2$. The clutter $\mathcal{H}_{12}$ is CA if and only if one of the following three cases holds.

\begin{enumerate}
    \item There is some intersection $A_i\cap A_{i+1}$ contained in both $B_1$ and $B_2$, and $A_i\cap B_1\subset A_i\cap B_2,\ A_{i+1}\cap B_2\subset A_{i+1}\cap B_1$.

    \item $B_1$ and $B_2$ contain $A_{i-1}\cap A_i$ and $A_i\cap A_{i+1}$ respectively. And, either $B_1$ and $B_2$ are disjoint, or the set $A_i\setminus (A_{i-1}\cup A_{i+1})$ is contained in the union $B_1\cup B_2$.

    \item $B_1$ and $B_2$ contain $A_{i-1}\cup A_i$ and $A_{i+1}\cap A_{i+2}$ respectively. And, either $B_1$ and $B_2$ are disjoint, or the set $A_i\cap A_{i+1}$ is contained in the union $B_1\cup B_2$.

\end{enumerate}

\begin{proof}

Let $A_i\cap A_{i+1}\subset B_1\subset A_i\cup A_{i+1}$ and $A_j\cap A_{j+1}\subset B_2\subset A_j\cup A_{j+1}$ by Lemma \ref{lem:intuni}. Let $T_{l,-}, T_{l,+}$ be the blocks that have nonempty intersections with $B_l$ at each end of $B_l$. An arc ordering of $\mathcal{H}_l$ is determined by orders of $T_{l,-}$, $T_{l,+}$ such that $B_l\cap T_{l,-}$ (resp. $B_l\cap T_{l,+}$) is an interval sharing the largest (resp. smallest) endpoint with $T_{l,-}$ (resp. $T_{l,+}$). The clutter $\mathcal{H}_{12}$ is CA if and only if there is a common arc ordering of $\mathcal{H}_1$ and $\mathcal{H}_2$. Therefore, no such common arc ordering exists if and only if 

\begin{enumerate}
    \item $T=T_{1,-}=T_{2,-}$, and $B_1\cap T, B_2\cap T$ are incomparable;

    \item $T=T_{1,+}=T_{2,-}$, and $B_1\cap B_2\cap T\neq\emptyset$, $T\neq (B_1\cap T)\cup(B_2\cap T)$.
\end{enumerate}

$T=T_{1,-}=T_{2,-}$ is possible only if $i=j$. So in this case a common interval ordering of $T$ exists if and only if $A_i\cap B_1\subset A_i\cap B_2$ (or $A_i\cap B_1\supset A_i\cap B_2$). Since $B_1, B_2$ are incomparable, this implies that $A_{i+1}\cap B_2\subset A_{i+1}\cap B_1$ (or $A_{i+1}\cap B_2\supset A_{i+1}\cap B_1$).

$T=T_{1,+}=T_{2,-}$ is possible only if $j=i+1,\ T=A_{i+1}\setminus(A_i\cup A_{i+2})$ or $j=i+2,\ T=A_{i+1}\cap A_{i+2}$.  This follows that case 2 or 3 must be true if there is a common arc ordering of $\mathcal{H}_1$ and $\mathcal{H}_2$.
    
\end{proof}
    
\end{lemma}

\begin{lemma}\label{lem:flag}

Let $V$ be a finite set of cardinality $n$. Suppose that there are two flags $A_1\subset A_2\subset\cdots\subset A_m$ and $B_1\supset B_2\supset\cdots\supset B_p$ of subsets of $V$. For any $i\in [m]$ and $j\in [p]$, assume that \begin{equation}
A_i\cap B_j\neq\emptyset \Rightarrow A_i\cup B_j=V. \label{eq:1}
\end{equation} Then, there is an interval ordering of the hypergraph $\mathcal{H}$ with the set of vertices $V$ and edges $E=\{A_1,\ldots,A_m,B_1,\ldots,B_p\}$ that maps each $A_i$ to the initial interval $[1,|A_i|]$ and each $B_j$ to the terminal interval $[n-|B_j|+1,n]$.

\begin{proof}

Suppose that $|A_{i+1}\setminus A_i|>1$. Choose $x\in A_{i+1}\setminus A_i$ that is contained in $B_j$ for $j=\max\{l\mid A_{i+1}\cap B_l\neq\emptyset\}$. By adjoining $A_{i+1}\setminus\{x\}$ to the original flag $A_1\subset\cdots\subset A_m$, the assumption (\ref{eq:1}) is preserved for this new flag. Therefore, it is possible to complete both flags while maintaining the property (\ref{eq:1}). Assuming both flags are complete, the bijection sending the unique element of $A_i \cap B_{n-i+1}$ to $i$ meets all requirements.
    
\end{proof}
    
\end{lemma}

Now, we prove the following theorem, which roughly states that ``local compatibility implies global compatibility for CA clutters''.

\begin{theorem}\label{thm:locglo}

Let $\mathcal{H}'=(V,E=\{A_1,\ldots, A_m, B_1,\ldots,B_p\})$ be a clutter with the partial clutter $\mathcal{H}=(V,\{A_1,\ldots,A_m\})$ that is a non-interval, CA clutter with every vertex contained in at most two edges. If the partial clutter $(V,\{A_1,\ldots,A_m,B_r,B_s\})$ is CA for all $r,s\in [p]$, then the clutter $\mathcal{H}$ is CA.

\begin{proof}

By Lemma \ref{lem:intuni}, the set $\{B_1,\ldots,B_p\}$ of edges can be partitioned into subsets $\mathcal{B}_i:=\{B_r\mid A_i\cap A_{i+1}\subset B_r\subset A_i\cup A_{i+1}\}$ for $i=1,\ldots,m$. Since the clutter $(V,\{A_1,\ldots,A_m,B_r,B_s\})$ is CA for every $B_r,B_s\in \mathcal{B}_i$, $A_i\cap B_r\subset A_i\cap B_s$ and $A_{i+1}\cap B_s\subset A_{i+1}\cap B_r$  by Lemma \ref{lem:incomp}. Therefore, the collection of intersections $\mathcal{F}_i(T):=\{B_r\cap T\mid B_r\in\mathcal{B}_i\}$ is a flag of $T$ for each block $T$.

Suppose that $T=A_i\cap A_{i+1}$. Then only nontrivial flags are $\mathcal{F}_{i-1}(T)$ and $\mathcal{F}_{i+1}(T)$. By the third case of Lemma \ref{lem:incomp}, these two flags satisfy the assumption of Lemma \ref{lem:flag}. By Lemma \ref{lem:flag}, there is an order of $T$ such that $\mathcal{F}_{i-1}(T)$ is a flag of initial intervals and $\mathcal{F}_{i+1}(T)$ is a flag of terminal intervals. If $T=A_i\setminus(A_{i-1}\cup A_{i+1})$, then only nontrivial flags are $\mathcal{F}_i(T)$ and $\mathcal{F}_{i+1}(T)$. Again by Lemma \ref{lem:flag}, there is an order of $T$ sending $\mathcal{F}_i(T)$ to a flag of initial intervals and $\mathcal{F}_{i+1}(T)$ to a flag of terminal intervals. Hence, the original hypergraph has an arc ordering obtained by cyclically concatenating the linear orders of each block.
    
\end{proof}
    
\end{theorem}

\subsection{Paving, non-positroid excluded minor of interval positroids with more than four edges}

First, we prove that Theorem \ref{thm:locglo} holds in a stronger sense for paving clutters that are excluded minors of interval positroids.

\begin{lemma}\label{lem:exclocglo}

Let $M$ be a paving excluded minor of interval positroids with a partial clutter $O_{a_1,\ldots,a_m;d}$. Suppose that for every edge $B$ of $M$, the clutter $O_{a_1,\ldots,a_m;d}+B$ is CA. Then, the clutter of $M$ is CA and therefore $M$ is a positroid.

\begin{proof}

By Lemma \ref{lem:incomp}, it is enough to prove that the clutter obtained by adjoining two edges to $O_{a_1,\ldots,a_m;d}$ is CA. Suppose that the clutter $O_{a_1,\ldots,a_m;d}+B_1+B_2$ is not CA for some $A_i\cap A_{i+1}\subset B_1,B_2\subset A_i\cup A_{i+1}$. Then, there is some block $T$ of $O_{a_1,\ldots,a_m;d}$ such that $T\cap B_1$ and $T\cap B_2$ are incomparable. Without losing generality, suppose that $T$ is contained in $A_i\setminus A_{i+1}$. The contraction $M/(A_i\cap A_{i+1})$ contains four edges $A_i'$, $A_{i+1}'$, $B_1'$, and $B_2'$ where $X':=X\setminus(A_i\cap A_{i+1})$. Since $B_1'$ and $B_2'$ are not contained in $A_i'$, both intersections $B_1'\cap A_{i+1}'$ and $B_2'\cap A_{i+1}'$ are nonempty. If $B_1'\cap B_2'\cap A_{i+1}'$ is nonempty, then the partial clutter $A_i'$, $B_1'$, and  $B_2'$ is a non-interval clutter of type 1 which is a contradiction. Thus, we may assume that $B_1'\cap A_{i+1}'$ and $B_2'\cap A_{i+1}'$ are disjoint. For the same reason, $B_1'\cap A_i'$ and $B_2'\cap A_i'$ must be disjoint. Thus, the four edges $B_1'$, $B_2'$, $A_i'$, and $A_{i+1}'$ form a non-interval clutter of type 5. Therefore, the contraction $M/(A_i\cap A_{i+1})$ is not an interval positroid which is a contradiction.

Now, suppose that the clutter $O_{a_1,\ldots,a_m;d}+B_1+B_2$ is not CA, and for some $i\neq j$, $A_i\cap A_{i+1}\subset B_1\subset A_i\cup A_{i+1}$ and $A_j\cap A_{j+1}\subset B_1\subset A_j\cup A_{j+1}$. Then, there is a block $T$ such that $B_1\cap B_2\cap T$ is nonempty and $T$ is not contained in $B_1\cup B_2$. Then, the clutter of the deletion $M-(T\setminus (B_1\cup B_2))$ is non-interval which is a contradiction. Therefore, $O_{a_1,\ldots,a_m;d}+B_1+B_2$ is CA for every $B_1, B_2$. This implies that the original matroid is a positroid.

\end{proof}
    
\end{lemma}

\begin{corollary}\label{cor}

Let $M$ be a paving excluded minor of interval positroids satisfying the following conditions:

\begin{enumerate}
    \item Edges of $M$ are $A_1,\ldots,A_m,B_1,\ldots,B_p$ for $p\ge 2$.

    \item The partial clutter consisting of edges $A_1,\ldots,A_m$ is $O_{a_1,\ldots,a_m;d}$.

    \item $M$ is not a positroid.
\end{enumerate}

Then, for some $r\in [p]$, the matroid obtained by relaxation of $B_r$ is a paving excluded minor of interval positroids that is not a positroid.

\begin{proof}

If the relaxation of $B_r$ is a positroid for all $r\in [p]$, then every partial clutter $A_1,\ldots,A_m,B_s$ is CA for all $s\in [p]$. Hence, the clutter of the original matroid $M$ is CA by Lemma \ref{lem:exclocglo}. Therefore, some relaxation of $B_r$ is a paving excluded minor of interval positroids that is not a positroid.
    
\end{proof}
    
\end{corollary}

Therefore, any paving, non-positroid excluded minor of interval positroids with five edges can be obtained either by adding one edge to $O_{a_1,\ldots,a_4;d}$ or by adding one edge to non-positroid paving excluded minors with four edges in FIGURE \ref{fig:ex4}. The next three propositions rule out three possibilities.

\begin{proposition}\label{prop:notO}

Suppose that $M$ is a non-positroid, paving excluded minor of interval positroids of rank $d$ with edges $A_1,\ldots,A_m,B$ for some $m\ge 4$. Then, the partial clutter obtained by removing one edge $B$ is not $O_{a_1,\ldots,a_m;d}$.

\begin{proof}

Suppose that the partial clutter obtained by removing $B$ is $O_{a_1,\ldots,a_m;d}$.

First, assume that $B$ is contained in $A_1\cup A_2$. If $A_1\cap A_2$ is not contained in $B$, then three edges $A_1,A_2,B$ form a non-interval CA clutter. Since any element $x\in A_3\cap A_4$ is not contained in $A_1,A_2,B$, the deletion $M-x$ is not an interval positroid. Hence, $A_1\cap A_2\subset B$. Since $M$ is not a positroid, we may assume that $B$ and $A_2\setminus (A_1\cup A_3)$ overlap and $B\cap A_2\cap A_3\neq \emptyset$. Then $A_2,A_3,B$ have a common point and each has exclusive elements which implies that $B,C,D_2$ is a partial clutter of type 2. Since $A_1$ is not contained in $B$ there is an element $x$ of $A_1$ not contained in $A_2\cup A_3\cup B$. Therefore, the deletion $M-x$ is not an interval positroid which is a contradiction.

Now, suppose that $B$ contains some element of $A_1\cap A_2$ and is not contained in $A_1\cup A_2$. Then, either three edges $A_1,A_2,B$ form a non-interval partial clutter of type 2, or we may assume that $A_2\setminus A_1$ is contained in $B$ by symmetry. If $A_1,A_2,B$ form a non-interval partial clutter of type 2, then $B$ must contain all elements not contained in $A_1\cup A_2$. Especially, $B$ contains $A_3\cap A_4$. And if $A_2\setminus A_1$ is contained in $B$, then $A_2\cap A_3$ is contained in $B$. Therefore, it can be assumed that $B$ contains some intersection $A_1\cap A_2$ and not contained in $A_1\cup A_2$ after reindexing if necessary. Thus, $B$ must contain all elements not contained in $A_1\cup A_2$. This implies that $A_3\cap A_4$ is contained in $B$ and therefore $B$ must contain all elements not contained in $A_3\cup A_4$. If $m\neq 4$, this implies that $A_1$ is contained in $B$ which is a contradiction. So, assume that $m=4$. Since $A_1,A_2$ are not contained in $B$, there are some element $x\in A_2\cap A_3$ and $y\in A_1\cap A_4$ not contained in $B$. This implies that $A_2,A_3,B$ is a non-interval, CA partial clutter. Thus, the deletion $M-y$ is not an interval positroid. This is a contradiction.

Therefore, $B$ must not contain an element contained in some intersection $A_i\cap A_{i+1}$. If $B$ contains elements of $A_i$ and $A_{i+1}$, then $A_i,A_{i+1},B$ form a non-interval, CA clutter. Then, the deletion $M-(A_{i+2}\cap A_{i+3})$ is not an interval positroid. Thus, we may assume that $B$ contains an element of $A_1$ and some other $A_i$ for $i\neq 2,m$. Let $\alpha$ be the smallest index such that $i\neq 1, A_i\cap B\neq\emptyset$ and $\beta$ be the largest index such that $A_i\cap B\neq\emptyset$. If $m> 4$, then either the non-interval clutter of type 5 $A_1,\ldots, A_\alpha, B$ or $A_1, A_m, A_{m-1},\ldots, A_\beta, B$ does not cover the vertex set $V$. This contradicts that $M$ is an excluded minor of interval positroids. Finally, suppose that $m=4$. Then, $A_1,A_2,A_3,B$ form a non-interval CA clutter. Hence, $A_4$ must be contained in $A_1\cup A_3$. By symmetry, $A_2$ is also contained in $A_1\cup A_3$. Therefore, $A_1\cup A_3$ contains three disjoint edges $A_2$, $A_4$, and $B$. This forces one of two edges $A_1,A_3$ to be large. But, any large edge of $O_{a_1,a_2,a_3,a_4;d}$ can not contain an element exclusively. Therefore, $B$ is contained in either $A_1$ or in $A_3$ which is a contradiction.

\end{proof}
    
\end{proposition}

\begin{proposition}

Suppose that $M$ is a non-positroid, paving excluded minor of interval positroids with five distinct edges $A,B,C,D_1,D_2$. Then, the partial clutter obtained by removing one edge $D_2$ is not $\Sh_{a,b,c;d}$.

\begin{proof}

Suppose that the partial clutter consisting of four edges $A,B,C,D_1$ is $\Sh_{a,b,c;d}$.

~
\paragraph{\textbf{Case 1}: $D_2$ contains some element $x\in D_1\setminus(A\cup B\cup C)$.}$ $\newline

Since $D_2$ is not contained in $D_1$, we may assume that $D_2 \cap (A \setminus D_1)$ is nonempty due to symmetry. Consequently, $A$, $D_1$, and $D_2$ form a triangle unless $A \cap D_1 \subset D_2$. Therefore, we have two possibilities: either $A \cap D_1 \subset D_2$, or $B \setminus D_1$ and $C \setminus D_1$ are subsets of $D_2$.

In the case where $A \cap D_1 \subset D_2$, there exists an element $y \in A \setminus D_1$ that is not contained in $D_2$, as $A$ is not entirely contained in $D_2$. Consequently, the deletion $M - y$ results in a clutter with four edges: $B$, $C$, $D_1$, and $D_2$. This clutter is non-interval by Lemma \ref{lem:tri}.

If $B \setminus D_1$ and $C \setminus D_1$ are subsets of $D_2$, then $B \cap D_1$ is not contained in $D_2$. This implies that the partial clutter $B, D_1, D_2$ forms a triangle. Consequently, $A \setminus D_1$ is contained in $D_2$. Thus, three disjoint small edges—$A$, $B$, and $C$—are contained in the union $D_1 \cup D_2$. This implies that either $D_1$ or $D_2$ is large. However, this leads to a contradiction, as the deletion $M - x$ would then have four edges forming a clutter of type 4, which is impossible.

~
\paragraph{\textbf{Case 2}: $D_2$ contains some element $x\in A\cap D_1$ and contained in $A\cup D_1$.}$ $\newline

Since $A, D_1, D_2$ do not cover the entire vertex set, they cannot form a triangle. Consequently, $A \cap D_1$ must be contained in $D_2$. Let $x$ be an element in $A$ that is not contained in $D_1$ because $A$ is not included in $D_2$. The deletion $M - x$ results in a set of four edges: $B, C, D_1, D_2$. Once again, this clutter is non-interval by Lemma \ref{lem:tri}.

~
\paragraph{\textbf{Case 3}: $D_2$ contains some element $x\in A\cap D_1$ and not contained in $A\cup D_1$.}$ $\newline

Let $y\in B\setminus D_1$ is contained in $D_2$. The partial clutter $B,D_1,D_2$ is non-interval unless $B\cap D_1$ is contained in $D_2$. Therefore, either $B\cap D_1$ or $(A\cup C)\setminus D_1$ is contained in $D_2$.

Suppose $B\cap D_1$ is contained in $D_2$. Then, there is an element $x$ of $B$ not contained in $D_2$. Previous cases imply that $D_2$ is not contained in $B\cup D_1$. Hence, $D_2$ contains an element in $(A\cup C)\setminus D_1$. Then, the partial clutter $B,D_1,D_2$ is a clutter of type 2. Therefore, $(A\cup C)\setminus D_1$ is contained in $D_2$.

Now, it remains to prove that $((A\cup C)\setminus D_1)\subset D_2$ induces a contradiction. The partial clutter with four edges $A,C,D_1,D_2$ is non-interval, as it has a subhypergraph that is a square. Therefore $B\setminus D_1$ is also contained in $D_2$. Since $A,B,C$ are not contained in $D_2$, $A\cap D_2, B\cap D_2, C\cap D_2$ are not entirely contained in $D_2$. Again, either $D_1$ or $D_2$ is large. If $D_1$ is large, at least one of the three intersections $A\cap D_1, B\cap D_1, C\cap D_1$ contains more than one element because $d\ge 3$. Suppose that $1<|A\cap D_1|$. The deletion $M-y$ by an element in $y\in (A\cap D_1)\setminus D_2$ results in a paving matroid with four edges $B,C,D_1\setminus\{y\},D_2$ which is not an interval positroid. By symmetry, the same contradiction is induced if $D_2$ is large.

~
\paragraph{\textbf{Case 4}: $D_1\cap D_2=\emptyset$.}$ $\newline

Combining three previous cases, the only remaining possibility is when $D_1$ and $D_2$ are disjoint. The edge $D_2$ must have nonempty intersection with at least two edges among $A,B,C$, as $D_2$ is not contained in another edge. Let us assume that $A\cap D_2$ and $B\cap D_2$ are nonempty. Then, the partial clutter obtained by removing $C$ is a non-interval CA clutter. Therefore, $C\setminus D_1$ is also contained in $D_2$. By the same argument, $A\setminus D_1$ and $B\setminus D_1$ are also contained in $D_2$. Therefore, $V=A\cup B\cup C=D_1\cup D_2$. This implies that either $D_1$ or $D_2$ is large. This induces a contradiction by the same argument as in case 3.

\end{proof}
    
\end{proposition}

\begin{proposition}\label{prop:notY}

Suppose that $M$ is a non-positroid, paving excluded minor of interval positroids with five distinct edges $A,B,C_1,C_2,D$. Then, the partial clutter obtained by removing one edge $D$ is not $Y'_{a,b;d}$.

\begin{proof}

Suppose that the partial clutter with three edges $A,B,C_1$ is $Y_{a,b;d}$ and $C_2=(A\setminus C_1)\sqcup(C_1\setminus(A\cup B))\sqcup(B\setminus C_1)$. Then, the partial clutter with four edges $A,B,C_1,D$ is also a paving, non-positroid excluded minor of interval positroids. By Theorem \ref{thm:Y4}, $D$ is either $(A\setminus C_1)\sqcup(C_1\setminus(A\cup B))\sqcup(B\setminus C_1)$ or $A,B,D$ form a clutter $Y_{a',b';d}$. Since $D$ is not equal to $C_2$, the latter case must hold. Now, the partial clutter with four edges $A,B,C_2,D$ is also a paving, non-positroid excluded minor of interval positroids. Since $C_2$ does not contain the common point, $C_2$ must be $(A\setminus D)\sqcup(D\setminus(A\cup B))\sqcup(B\setminus D)$ which implies that $D=C_1$. This is a contradiction.
    
\end{proof}
    
\end{proposition}

Next theorem characterizes all non-positroid, paving excluded minor of interval positroids obtained by adding edges to $Y_{a,b;d}$.

\begin{theorem}

Suppose that $M$ is a paving excluded minor of interval positroids of rank $d$ with edges $A,B,C_1,C_2,\ldots,C_m$ for some $m\ge 3$. The partial clutter obtained by removing edges $C_2,\ldots,C_m$ is $Y_{a_1,b_1;d}$ with the common point $v$ if and only if following conditions hold:

\begin{enumerate}

    \item The partial clutter $A,B,C_i$ is the clutter $Y_{a_i,b_i;d}$ for each $i=1,\ldots,m$.
    
    \item The clutter $(V\setminus\{v\},\{A\setminus\{v\},B\setminus\{v\},C_1\setminus\{v\},\ldots,C_m\setminus\{v\}\})$ is an interval $(d-1)$-paving clutter.

\end{enumerate}

\begin{proof}

The proof of sufficiency is omitted as it can be checked routinely.

The partial clutter with four edges $A,B,C_1,C_i$ must be $Y_{a_1,a_i,b_1,b_i;d}$ for $i=2,\ldots,m$ by Theorem \ref{thm:Y4} and Proposition \ref{prop:notY}. Hence, the first condition holds. The second condition follows from the fact that the contraction $M/v$ is an interval positroid of rank $d-1$.
    
\end{proof}
    
\end{theorem}

\begin{notation}

The excluded minor above will be denoted by $Y_{a_1,\ldots,a_m,b_1,\ldots,b_m;d}$ where $a_i=|A\cap C_i|-1, b_i=|B\cap C_i|-1$.

\end{notation}

Now, we classify paving, non-positroid, $O_{a_1,a_2,a_3;d}$-base excluded minors of interval positroids. Recall the notation \begin{align*}
    S&=A\setminus(B\cup C), S'=B\cap C, \\
    T&=B\setminus(C\cup A), T'=C\cap A, \\
    U&=C\setminus(A\cup B), U'=A\cap B, \\
    a&=|S'|, b=|T'|, c=|U'|.
\end{align*}

\begin{theorem}\label{thm:delone}

Suppose that $M$ is a non-positroid, paving excluded minor of interval positroids with five distinct edges $A,B,C,D_1,D_2$. If the partial clutter obtained by removing one edge $D_2$ is $\Delta_{a,b,c}^1$, then $D_2$ must be $S'\cup T'\cup U'$. This paving excluded minor of interval positroids will be denoted by $\Delta_{a,b,c}^3$.

\begin{proof}

If the partial clutter $A,B,C,D_2$ is not a CA clutter, then $D_2$ must be either $S'\cup T'\cup U'$ or $S\cup T\cup U$ by Theorem \ref{thm:O4}. Since $D_2$ is not $D_1$, it must be $S'\cup T'\cup U'$.

Now, assume that the partial clutter with four edges $A,B,C,D_2$ is CA. By symmetry and Lemma \ref{lem:intuni}, we may assume the inclusion relations $S'\subset D_2\subset B\cup C$. If $D_2$ is the complement of the edge $A$, the matroid is $\Delta_{a,b,c}^3$. f $D_2$ is not the complement of the edge $A$, we can assume that the edge $D_2$ contains some elements $x \in T$ and $y \in T'$. Then $D_2$ contains two consecutive blocks $S'$ and $U$. Then, the partial clutter with three edges $B,D_1,D_2$ is a clutter of type 2. Hence, $D_2$ must contain $T'=A\cap C$. This implies that $C\subset D_2$ which is a contradiction.
    
\end{proof}

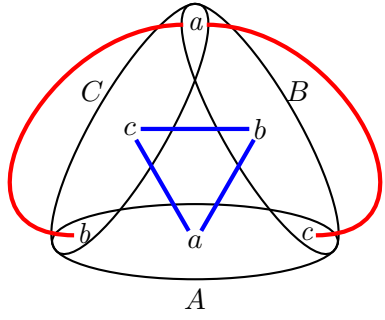
\begin{figure}
    \centering
    \begin{tikzpicture}
  \begin{scope}[rotate=300]
    \draw[thick] (30:1) ellipse [x radius=1.9, y radius=0.5] node[] (b) {};
  \end{scope}
  \begin{scope}[rotate=60]
    \draw[thick] (150:1) ellipse [x radius=1.9, y radius=0.5] node[] (c) {};
  \end{scope}
  \draw[thick] (270:2) ellipse [x radius=1.9, y radius=0.5] node[] (a) {};
  \node at (a) [] {$a$};
  \node at (a) [below, yshift=-0.5cm] {$A$};
  \node at (b) [] {$b$};
  \node at (b) [xshift=0.5cm, yshift=0.5cm] {$B$};
  \node at (c) [] {$c$};
  \node at (c) [xshift=-0.5cm, yshift=0.5cm] {$C$};
  \node at (90:0.9) {a};
  \node at (232:2.4) {b};
  \node at (308:2.4) {c};
  \draw[ultra thick, draw=red] (100:0.9) to[out=180, in=180, looseness=1.5] (230:2.5);
  \draw[ultra thick, draw=red] (80:0.9) to[out=0, in=0, looseness=1.5] (310:2.5);
  \draw[ultra thick, draw=blue, rounded corners=8pt] (a) -- (b) -- (c) -- (a);
\end{tikzpicture}
    \caption{$\Delta_{a,b,c;d}^3$}
    \label{fig:O5}
\end{figure}

\end{theorem}

Similarly, $\Delta^3_{a,b,c}$ is the only possible case if some relaxation is $\Delta_{a,b,c;k}^2$.

\begin{theorem}\label{thm:deltwo}

Suppose that $M$ is a non-positroid, paving excluded minor of interval positroids with five distinct edges $A,B,C,D_1,D_2$. If the partial clutter obtained by removing one edge $D_2$ is $\Delta_{a,b,c;k}^2$, then $k=0$ and $D_2$ must be $S\cup T\cup U$. Therefore, $M$ is $\Delta^3_{a,b,c}$

\begin{proof}

If the partial clutter $A,B,C,D_2$ is not a CA clutter, then $D_2$ must be either $S'\cup T'\cup U'$ or $S\cup T\cup U$ (assuming $k=0$) by Theorem \ref{thm:O4}. Since $D_2$ is not equal to $D_1$, $D_2$ must be $S\cup T\cup U$ and $k=0$.

Now, assume that the partial clutter with four edges $A,B,C,D_2$ is CA. Assume that $S'\subset D_2\subset B\cup C$. Then $B,C,D_1,D_2$ is a non-interval clutter of type 3. Since $S$ is nonempty, the deletion $M-S$ is a proper minor of $M$ that is not an interval positroid.
    
\end{proof}
    
\end{theorem}

Finally, we end this section by showing that $Y_{a_1,\ldots,a_m,b_1,\ldots,b_m;d}, m\ge3$ are only possible paving, non-positroid excluded minors of interval positroids with more than five edges. By Proposition \ref{prop:notO}, no such excluded minor is obtained by adding an edge to $O_{a_1,\ldots,a_m;d}$. Moreover, Corollary \ref{cor} implies that some relaxation of such excluded minor remains non-positroid. Therefore, it is enough to prove that no new excluded minor is obtained by adding an edge to $\Delta_{a,b,c}^3$, as $\Delta_{a,b,c}^3$ and $Y_{a_1,a_2,a_3,b_1,b_2,b_3;d}$ are only non-positroid paving excluded minor of interval positroids with five edges.

\begin{proposition}

Let $M$ be a non-positroid, paving excluded minors of interval positroids with six distinct edges $A,B,C,D_1,D_2,D_3$. Then, the partial clutter $A,B,C,D_1,D_2$ is not $\Delta_{a,b,c}^3$.

\begin{proof}

Suppose that the partial clutter $A,B,C,D_1,D_2$ is $\Delta_{a,b,c}^3$ with $D_1=S\cup T\cup U$ and $D_2=S'\cup T'\cup U'$. By Theorem \ref{thm:delone}, the partial clutter $A,B,C,D_1,D_3$, which is also a paving, non-positroid excluded minor of interval positroids, must be $\Delta_{a,b,c}^3$. Thus $D_3$ must be $S'\cup T'\cup U'$ which is a contradiction since $D_2$ is not equal to $D_3$.
    
\end{proof}
    
\end{proposition}

\section{Summary and further remarks}\label{sec:sum}

To summarize, we list all paving, non-positroid excluded minors of interval positroids in Table \ref{tab}.

\begin{table}
    \centering
    \begin{tabular}{|P{6cm}|P{6cm}|}
    \hline
        Number of edges or dependent hyperplanes & non-positroid paving excluded minors of interval positroids\\
        \hline
        3 & $Y_{a,b;d}$\\
        \hline
        4 & $Y'_{a,b;d}, Y_{a_1,a_2,b_1,b_2;d}, \Delta_{a,b,c}^1, \Delta_{a,b,c;k}^2$\\
        \hline
        5 & $Y_{a_1,a_2,a_3,b_1,b_2,b_3;d}, \Delta_{a,b,c}^3$\\
        \hline
        6 & $Y_{a_1,\ldots,a_m,b_1,\ldots,b_m;d}$ for $m\ge 4$\\
        \hline
    \end{tabular}
    \caption{Non-positroid paving excluded minors of interval positroids}
    \label{tab}
\end{table}

A paving, non-positroid excluded minors of interval positroids is also an excluded minor of positroids. Moreover, while the dual of these excluded minors is not necessarily an excluded minor of interval positroids, it remains an excluded minor of positroids. The excluded minor $\Delta_{a,b,c}^1$, which is self-dual, is a special case of the example found by \citet[Example~4,5]{bonin2023characterization}. The dual of the excluded minor $Y_{a,b;d}$, when the edge $C$ is small, is a special case of the three-circuit case of the example \cite[Example~9]{bonin2023characterization}. The excluded minors of positroids $Y_{a_1,\ldots,a_m,b_1,\ldots,b_m;d}$ (for $m\ge 2$ or with some large $C_i$), $Y'_{a,b;d}, \Delta_{a,b,c;k}^2, \Delta_{a,b,c}^3$ and their duals are new excluded minors of positroids that have not appeared in the literature.

\bibliographystyle{unsrtnatneat}
\bibliography{ref}

\end{document}